\documentclass[a4paper,11pt]{amsart}
\usepackage{amsmath,amsthm,amssymb,hyperref,enumitem}
\usepackage[T1]{fontenc}
\usepackage{lmodern}

\newtheorem{theorem}{Theorem}[section]
\newtheorem{proposition}[theorem]{Proposition}
\newtheorem{lemma}[theorem]{Lemma}
\newtheorem{corollary}[theorem]{Corollary}

\theoremstyle{definition}
\newtheorem{remark}[theorem]{Remark}

\begin{document}

\title[Holomorphic interpolation and the Levi problem]{Holomorphic interpolation\\and the Levi problem}
\author{Yuta Kusakabe}
\address{Faculty of Mathematics, Kyushu University, Fukuoka 819-0395, Japan}
\email{kusakabe@math.kyushu-u.ac.jp}
\author{Takeo Ohsawa}
\address{Graduate School of Mathematics, Nagoya University, Nagoya 464-8602, Japan}
\email{ohsawa@math.nagoya-u.ac.jp}

\begin{abstract}
On a nonsingular locally Stein Riemann domain over a Stein space, we prove that holomorphic forms of top degree interpolate jets of any fixed order on every closed discrete set finite over the base. For unramified domains, this yields jet interpolation by holomorphic functions, answering a question of Col\c{t}oiu and Joi\c{t}a for nonsingular sources. More generally, we establish a twisted interpolation theorem over weakly pseudoconvex spaces admitting injective holomorphic immersions into projective space. We also give negative answers to the Levi problem over singular Stein spaces: the C\oe ur\'e--Loeb space yields a non-Stein locally Stein unramified Riemann domain over a normal Stein threefold, and a related construction gives a relatively compact, non-Stein locally Stein open subset of a normal Stein fourfold. Finally, the punctured anticanonical bundle of Forn\ae ss' ramified surface is a holomorphically separable increasing union of Stein manifolds that is not Stein, answering a question of Vigna Suria.
\end{abstract}

\maketitle

\section{Introduction}\label{sec:intro}

A holomorphic map $\pi\colon X\to Y$ is \emph{locally Stein} if every point of $Y$ has an open neighborhood $U$ such that $\pi^{-1}(U)$ is Stein. An open subset of $Y$ is locally Stein if its inclusion map is. A \emph{Riemann domain} is a holomorphic map with discrete fibers between complex spaces of the same pure dimension, with connected source; it is \emph{unramified} if it is a local biholomorphism, \emph{ramified} otherwise, and \emph{nonsingular} if its source is a manifold. Complex spaces are reduced, Hausdorff, second countable and finite-dimensional.

The Levi problem asks whether local Steinness implies Steinness. For unramified domains over $\mathbb{C}^n$ it was solved by Oka~\cite{Oka53}, and over Stein manifolds by Docquier and Grauert~\cite{DG}. Over $\mathbb{P}^n$, Fujita and Takeuchi proved that a locally Stein unramified domain is Stein unless its projection is a biholomorphism~\cite{Fuj63,Tak64}. The implication fails for ramified domains: Grauert gave a counterexample over $\mathbb{P}^2$~\cite[p.~377]{Gr63}, and Forn\ae ss constructed a locally Stein ramified surface $\pi\colon M\to\Delta\times\mathbb{C}$, where $\Delta$ is the unit disk, that is not Stein~\cite{F78}. For open subsets and unramified Riemann domains over singular Stein spaces, the problem has remained open~\cite{Si,C09}. Since a complex manifold is Stein if and only if holomorphic functions interpolate arbitrary values on every closed discrete set, we ask how much interpolation local Steinness still gives.

Functions on a complex space are smooth, or plurisubharmonic, if they are locally restrictions of such functions under local embeddings~\cite{FN80}, and a complex space is \emph{weakly pseudoconvex}, or weakly $1$-complete, if it carries a smooth plurisubharmonic exhaustion function. Let $Y$ be a weakly pseudoconvex space of pure dimension $n\ge1$ with an injective holomorphic immersion $\iota\colon Y\to\mathbb{P}^N$, that is, an injective holomorphic map that maps a neighborhood of every point isomorphically onto a closed subspace of an open subset of $\mathbb{P}^N$, not necessarily a homeomorphism onto its image, and put $\mathcal{O}_Y(1)=\iota^*\mathcal{O}_{\mathbb{P}^N}(1)$ and $\mathcal{O}_Y(m)=\mathcal{O}_Y(1)^{\otimes m}$. Examples are projective varieties, closed analytic subsets of $\mathbb{C}^N$ and, by the embedding theorems of Kodaira and of Takayama~\cite[Theorem~1.2]{Taka98}, weakly pseudoconvex manifolds with a positive line bundle.

Let $\pi\colon X\to Y$ be a nonsingular Riemann domain. A closed discrete set $S\subset X$ is \emph{finite over $Y$} if $S\cap\pi^{-1}(B)$ is finite for every compact set $B\subset Y$. Two points of a complex space are \emph{indistinguishable} if every smooth plurisubharmonic function takes the same value at them; for instance, the classes of indistinguishable points of $\mathbb{P}^1\times\mathbb{C}$ are the curves $\mathbb{P}^1\times\{t\}$, and those of a Stein space are its points. A line bundle and its sheaf of holomorphic sections share their notation, a subscript $x$ denotes the stalk, $\mathfrak{m}_x$ is the maximal ideal of $\mathcal{O}_{X,x}$, and the \emph{$k$-jet} $j^k_xs$ of a section $s$ of a line bundle $L$ at $x$ is its class in the space $L_x/\mathfrak{m}_x^{k+1}L_x$ of $k$-jets. On Forn\ae ss' surface, holomorphic functions do not even separate points (Theorem~\ref{thm:fornaess}), but holomorphic $2$-forms interpolate jets of any fixed order on closed discrete sets finite over the base.

\begin{theorem}\label{thm:main}
Suppose that $\pi$ is locally Stein. Let $L$ be a holomorphic line bundle on $X$ with a smooth hermitian metric of semipositive curvature, and let $k\ge0$, $p\ge1$ and $m\ge p(n+k)+1$ be integers. Let $S\subset X$ be a closed discrete set finite over $Y$ such that every class of indistinguishable points of $Y$ contains at most $p$ points of $\pi(S)$, and for $x\in S$ let $a_x$ be a $k$-jet at $x$ of a section of $K_X\otimes L\otimes\pi^*\mathcal{O}_Y(m)$. Then there is $s\in H^0(X,K_X\otimes L\otimes\pi^*\mathcal{O}_Y(m))$ such that
\begin{equation}\label{eq:interp}
j^k_xs=a_x\quad(x\in S),\qquad j^k_zs=0\quad\bigl(z\in\pi^{-1}(\pi(S))\setminus S\bigr).
\end{equation}
\end{theorem}

The bound on $m$ depends neither on $Y$, $\iota$ and the number of sheets, nor on $L$, and it is sharp for $p=1$: if $X=Y=\mathbb{P}^n$ and $L$ is trivial, the sections of $K_X\otimes\mathcal{O}_Y(m)=\mathcal{O}_{\mathbb{P}^n}(m-n-1)$ interpolate $k$-jets at one point only if $m\ge n+k+1$. The bounds on the classes and on the jet orders cannot be dropped: if $X=Y=\mathbb{P}^1\times\mathbb{C}$, embedded in $\mathbb{P}^3$ by the Segre map, a line bundle has the same degree $d$ on all curves $\mathbb{P}^1\times\{t\}$, and the restriction of a section to such a curve vanishes identically or has at most $d$ zeros, counted with multiplicity. Hence the sections neither interpolate arbitrary values on $\max(d+2,1)$ points of one such curve nor interpolate jets on $\{b\}\times\mathbb{N}$ whose orders tend to infinity.

The proof solves a $\bar\partial$-equation on the unramified part over a dense Stein open subset of the regular part of the complement of a hyperplane section, which is Stein by Docquier and Grauert, with weights pulled back from the base: linear forms provide the poles, and level sets of a plurisubharmonic exhaustion function of $Y$ localize them. Thus no plurisubharmonic exhaustion along the fibers is needed; its lack was pointed out by Col\c{t}oiu and Joi\c{t}a as an obstacle to $L^2$ methods on infinitely sheeted domains~\cite{CoJo20}. Similar $L^2$ arguments were used for locally Stein open subsets of Stein spaces in~\cite[Theorem~3.2.1]{FN80} and~\cite{Vaj19}.

If $Y$ is instead a Stein space, a proper injective holomorphic map $F\colon Y\to\mathbb{C}^{2n+1}$~\cite{N60} turns a locally Stein $\pi$ into a locally Stein Riemann domain $F\circ\pi$ over the closed analytic subset $F(Y)$ of $\mathbb{C}^{2n+1}\subset\mathbb{P}^{2n+1}$, on which $\mathcal{O}(1)$ is trivial and the classes of indistinguishable points are points. Hence Theorem~\ref{thm:main} gives $s\in H^0(X,K_X\otimes L)$ satisfying~\eqref{eq:interp} for every closed discrete set $S$ finite over $Y$ and every $k\ge0$ (Section~\ref{sec:consequences}); we call this the Stein case of Theorem~\ref{thm:main}.

\begin{corollary}\label{cor:main}
Let $\pi\colon X\to Y$ be a nonsingular locally Stein Riemann domain.
\begin{enumerate}[label=\textup{(\roman*)}]
\item If $Y$ is as in Theorem~\textup{\ref{thm:main}}, then $K_X\otimes\pi^*\mathcal{O}_Y(n+1)$ is generated by global sections. For every $m\ge n+2$, sections of $K_X\otimes\pi^*\mathcal{O}_Y(m)$ give a holomorphic map $\Phi\colon X\to\mathbb{P}^{n+1}$ such that $(\pi,\Phi)\colon X\to Y\times\mathbb{P}^{n+1}$ is an injective holomorphic immersion. If $Y$ is a Stein space, the same holds for $K_X$ itself.
\item If $Y$ is a Stein space, and if $\pi$ is unramified or $K_X^{-1}$ carries a smooth hermitian metric of semipositive curvature, then holomorphic functions on $X$ interpolate jets of any fixed order on closed discrete sets finite over $Y$.
\end{enumerate}
\end{corollary}

Part~(i) is a relative analog of Takayama's theorem that, for a positive line bundle $L$ on an $n$-dimensional weakly pseudoconvex manifold $Z$, the adjoint bundles $K_Z\otimes L^{\otimes m}$ are ample for all $m>n(n+1)/2$~\cite[Theorem~1.2]{Taka98}, and of its extension by Watanabe to the regular locus of a weakly pseudoconvex space~\cite[Theorem~1.1]{Wat}; here the base supplies the positivity and $\pi$ the remaining coordinates. Part~(ii) answers affirmatively, for nonsingular sources, a question of Col\c{t}oiu and Joi\c{t}a, who asked whether holomorphic functions separate the points of every locally Stein unramified Riemann domain over a Stein space~\cite[Problem~2]{CoJo20}. Following~\cite[Remark~2]{CoJo20}, it also answers their Problem~3 for nonsingular sources.

\begin{corollary}\label{cor:envelope}
Let $\pi\colon X\to Y$ be a nonsingular locally Stein unramified Riemann domain over a Stein space. Then $X$ is Stein if and only if $\pi$ has an envelope of holomorphy, that is, a Stein unramified Riemann domain $p\colon\tilde X\to Y$ with a holomorphic map $\tau\colon X\to\tilde X$ such that $p\circ\tau=\pi$ and every $f\in\mathcal{O}(X)$ is $\tilde f\circ\tau$ for a unique $\tilde f\in\mathcal{O}(\tilde X)$.
\end{corollary}

We follow~\cite[Theorem~3.3]{CD00}, replacing~\cite[Theorem~3.1.1]{FN80} by Preda's theorem~\cite{Pre15}. Since the ambient space is the envelope itself, it may be normalized, and no relative compactness is needed.

Over Stein spaces, interpolation on all closed discrete sets detects Steinness even for sections of a line bundle.

\begin{proposition}\label{prop:steinness}
Let $\pi\colon X\to Y$ be a nonsingular Riemann domain over a Stein space. Then $X$ is Stein if and only if the sections of some holomorphic line bundle on $X$ interpolate arbitrary values on every closed discrete subset of $X$.
\end{proposition}

Neither local Steinness nor positivity is assumed. Thus a nonsingular locally Stein Riemann domain over a Stein space is Stein precisely when holomorphic $n$-forms interpolate arbitrary values on every closed discrete subset, not only on those finite over $Y$.

We apply these results to Forn\ae ss' ramified surface and the C\oe ur\'e--Loeb space.

\begin{theorem}\label{thm:fornaess}
Let $\pi\colon M\to\Delta\times\mathbb{C}$ be Forn\ae ss' ramified surface.
\begin{enumerate}[label=\textup{(\alph*)}]
\item Every holomorphic function on $M$ factors through $\pi$. In particular, $M$ is not holomorphically separable.
\item The punctured anticanonical bundle $(K_M^{-1})^\times$ is a holomorphically separable increasing union of Stein open subsets that is not Stein.
\end{enumerate}
\end{theorem}

In contrast to holomorphic functions, holomorphic $2$-forms on $M$ separate points by the Stein case of Theorem~\ref{thm:main}.

Part~(b) answers a question of Vigna Suria~\cite{VS84}, who showed, using~\cite{F77}, that an increasing union of Stein manifolds need not be holomorphically separable, and asked whether a holomorphically separable such union must be Stein. Moreover, holomorphic functions and vector fields on $(K_M^{-1})^\times$ have arbitrary jets at any finite set, and, as an increasing union of Stein manifolds, it satisfies the Kontinuit\"atssatz for continuous families of disks (Proposition~\ref{prop:siu}); thus it also answers negatively the weaker form of a question of Siu~\cite[Question~(8.5)]{Si}. By Watanabe's solution of the union problem for weakly pseudoconvex manifolds~\cite[Theorem~1.7]{Wat}, $(K_M^{-1})^\times$ is not weakly pseudoconvex. We do not know whether $(K_M^{-1})^\times$ is biholomorphic to an open subset of a Stein space. A positive answer would give a counterexample to the union problem for Stein spaces~\cite[Question~(5.5)]{Si}, \cite[Section~1.2]{FN80}, \cite[Question~2]{C09}; a negative answer would give a nonsingular counterexample to the open immersion problem, which is open for normal spaces~\cite{CoJo13,Pre26}.

The C\oe ur\'e--Loeb space $E$ is a non-Stein holomorphic fiber bundle over $\mathbb{C}^*$ with fiber a bounded Stein domain in $\mathbb{C}^2$, a counterexample to the Serre problem~\cite{CL}. Col\c{t}oiu and Diederich observed, using results of Siu~\cite{Siu76b}, that $E$ is a locally Stein ramified Riemann domain over $\mathbb{C}^*\times\mathbb{C}^2$ on which holomorphic functions separate points and give local coordinates~\cite[Section~4, Remark~3]{CD07}. Since $K_E$ is trivial (Proposition~\ref{prop:CL-canonical}), Corollary~\ref{cor:main}(ii) strengthens this conclusion to interpolation of jets on closed discrete sets finite over $\mathbb{C}^*\times\mathbb{C}^2$, although $E$ is not Stein.

The space $E$ also bears on the Levi problem over singular Stein spaces, which asks, for a Stein space $Y$, whether \textup{(1)}~every locally Stein open subset of $Y$ and \textup{(2)}~every locally Stein unramified Riemann domain over $Y$ is Stein~\cite[Problem~(4.2)]{Si}, \cite[Questions~1 and~2]{CD07}, \cite[Questions~1 and~1$'$]{C09}. Both have affirmative answers for isolated singularities, by Andreotti and Narasimhan~\cite{AN} and by Col\c{t}oiu and Diederich~\cite{CD07}, and in dimension two~\cite[p.~362]{AN}, \cite{CoJo20}. Col\c{t}oiu proposed to derive counterexamples to~(1) from counterexamples to the Serre problem~\cite[Section~6]{C09}.

\begin{theorem}\label{thm:CL}
There is a normal Stein surface $Y_A$ with one cusp singularity such that the C\oe ur\'e--Loeb space $E$ is a locally Stein unramified Riemann domain over $\mathbb{C}^*\times Y_A$, namely an infinite cyclic covering of $\mathbb{C}^*\times Y_{A,\mathrm{reg}}$. In particular, a locally Stein unramified Riemann domain over a Stein space of dimension three need not be Stein.
\end{theorem}

As $E$ is unramified over $\mathbb{C}^*\times Y_A$, Corollary~\ref{cor:main}(ii) shows that holomorphic functions on $E$ still interpolate jets on closed discrete sets finite over $\mathbb{C}^*\times Y_A$. Answering a question of Col\c{t}oiu~\cite[Problem~7]{C19}, Preda proved that $E$ is an open subset of a Stein space~\cite[Theorem~3.2]{Pre26}. We refine the ambient space to a normal Stein threefold with one singular point (Theorem~\ref{thm:CLembed}), and thickening the graph of the bundle projection answers~(1) negatively in dimension four.

\begin{theorem}\label{thm:open}
There exist a normal Stein threefold $Z$ with exactly one singular point and a connected, locally Stein open subset $W\Subset\mathbb{C}^*\times Z$ that is contained in $(\mathbb{C}^*\times Z)_{\mathrm{reg}}$ and is not Stein.
\end{theorem}

Taking products with polydisks, (2) and (1) have negative answers in every dimension at least three and four, respectively. Thus (2) is settled in all dimensions, while (1) remains open only in dimension three, where Col\c{t}oiu conjectured an affirmative answer~\cite[Conjecture~1]{C09}. If this conjecture holds, then $W\cap H$ is Stein for every hypersurface $H$ of $\mathbb{C}^*\times Z$, as a locally Stein open subset of the three-dimensional Stein space $H$, so $W$ is a counterexample to the hypersurface section problem in dimension four~\cite[Section~5]{C09}. Theorem~\ref{thm:open} also answers negatively the versions of~(1) for relatively compact open subsets~\cite[Section~1.1]{FN80} and for products of a Stein space with an isolated singularity and a Stein manifold~\cite[Section~6]{C09}; see also Remark~\ref{rem:gluing}.

\section{Preliminaries}\label{sec:prelim}

\subsection{Riemann domains}
Let $\pi\colon X\to Y$ be a holomorphic map with discrete fibers. If a closed discrete set $S\subset X$ is finite over $Y$, then $\pi(S)$ is closed and discrete as well, since $\pi(S)\cap B=\pi(S\cap\pi^{-1}(B))$ is finite for every compact set $B\subset Y$.

\begin{lemma}\label{lem:finrep}
Let $\pi\colon X\to Y$ be a continuous map with discrete fibers between locally compact Hausdorff spaces, let $x\in X$, $y=\pi(x)$, and let $W'$ and $P'$ be neighborhoods of $x$ and $y$. There are a relatively compact open neighborhood $W\subset W'$ of $x$ and an open neighborhood $P\subset P'$ of $y$ such that $\overline W\cap\pi^{-1}(y)=\{x\}$, $\pi(W)\subset P$ and $\pi|_W\colon W\to P$ is proper. For every open neighborhood $P_1\subset P$ of $y$, the restriction $W\cap\pi^{-1}(P_1)\to P_1$ is again proper.
\end{lemma}

\begin{proof}
Since $\pi^{-1}(y)$ is discrete, there is an open set $W_0\Subset W'$ containing $x$ with $\overline{W_0}\cap\pi^{-1}(y)=\{x\}$. The set $\pi(\partial W_0)$ is compact and does not contain $y$. Let $P\subset P'\setminus\pi(\partial W_0)$ be an open neighborhood of $y$ and $W=W_0\cap\pi^{-1}(P)$. For compact $Q\subset P$ we have $(\pi|_W)^{-1}(Q)=\overline{W_0}\cap\pi^{-1}(Q)$, because $\partial W_0\cap\pi^{-1}(Q)=\emptyset$, and this set is compact. The last assertion follows by restriction.
\end{proof}

We call $\pi|_W\colon W\to P$ a \emph{finite representative} of $\pi$ at $x$. If $\theta$ is a smooth function with compact support in $P$, then $(\theta\circ\pi)|_W$, extended by zero, is a smooth function on $X$ with compact support, since $W\cap\pi^{-1}(\operatorname{supp}\theta)$ is compact.

\subsection{\texorpdfstring{$L^2$ estimates for $n$-forms}{L2 estimates for n-forms}}
For a holomorphic line bundle $\mathcal{L}$ with a hermitian metric $h$ and an $\mathcal{L}$-valued $(n,0)$-form $u$, we write
\[
|u|^2_h=i^{n^2}h(u\wedge\bar u),
\]
a density that depends on $h$ but on no metric of the manifold. In local coordinates $w$ and a local frame $e$, if $u=f\,dw_1\wedge\dots\wedge dw_n\otimes e$, then $|u|^2_h$ is $|f|^2$ times a positive multiple of Lebesgue measure, the multiple being smooth if $h$ is. If $\omega$ is a K\"ahler form, $\alpha$ a $(0,1)$-form and $A$ an $\mathcal{L}$-valued $(n,0)$-form, then
\begin{equation}\label{eq:product}
|\alpha\wedge A|^2_{\omega,h}\,dV_\omega=|\alpha|^2_\omega\,|A|^2_h,
\end{equation}
as one checks in coordinates orthonormal for $\omega$ at a point. The following is the case $q=1$ of Demailly's $L^2$ estimate for $(n,q)$-forms on weakly pseudoconvex K\"ahler manifolds, in which the K\"ahler form need not be complete~\cite[Chapter~VIII, Theorem~6.5]{Dem}.

\begin{lemma}\label{lem:L2}
Let $\Omega$ be an $n$-dimensional weakly pseudoconvex manifold with a K\"ahler form $\omega$, let $\mathcal{L}$ be a holomorphic line bundle on $\Omega$ with a smooth hermitian metric $h$, and let $\Psi$ be a smooth function on $\Omega$ with $i\Theta(h)+i\partial\bar\partial\Psi\ge\omega$. For every $\bar\partial$-closed $\mathcal{L}$-valued $(n,1)$-form $v$ on $\Omega$ for which the right-hand side below is finite, there is an $\mathcal{L}$-valued $(n,0)$-form $u$ on $\Omega$ with $\bar\partial u=v$ and
\[
\int_\Omega|u|^2_he^{-\Psi}\le\int_\Omega|v|^2_{\omega,h}\,dV_\omega\,e^{-\Psi}.
\]
\end{lemma}

A function, or a section of a line bundle, that is holomorphic off a nowhere dense analytic subset and square integrable near it extends holomorphically across it, since it is $\bar\partial$-closed across the subset in the sense of distributions~\cite[Chapter~VIII, Lemma~7.3]{Dem}. Weights with a pole at a point force jets to vanish:

\begin{lemma}\label{lem:jets}
Let $n\ge1$, let $f$ be holomorphic on $B(0,\delta)\subset\mathbb{C}^n$, and let $\mu\ge n-1$ be an integer with $\int_{B(0,\delta)}|f|^2|w|^{-2\mu}\,d\lambda<\infty$. Then $f\in\mathfrak{m}_0^{\mu-n+1}$.
\end{lemma}

\begin{proof}
If $f\not\equiv0$, let $m=\operatorname{ord}_0f$ and let $P_m$ be its first nonzero homogeneous term. As $r\downarrow0$, the functions $r^{-m}f(r\xi)$ converge uniformly to $P_m(\xi)$ on the unit sphere, so
\[
\int_{|\xi|=1}|f(r\xi)|^2\,d\sigma\sim c r^{2m},\qquad c=\int_{|\xi|=1}|P_m(\xi)|^2\,d\sigma>0.
\]
Polar coordinates therefore give $2m-2\mu+2n>0$, or $m\ge\mu-n+1$.
\end{proof}

\subsection{Plurisubharmonic exhaustion functions}
By definition, indistinguishable points lie on the same level set of every smooth plurisubharmonic function. The following lemma separates all other points of a closed discrete set simultaneously by the level sets of one exhaustion function.

\begin{lemma}\label{lem:levels}
Let $Z$ be a complex space with a smooth plurisubharmonic exhaustion function $\varphi$, let $D\subset Z$ be a closed discrete set, and let $p\ge1$ be an integer such that every class of indistinguishable points of $Z$ contains at most $p$ points of $D$. Then $Z$ has a continuous plurisubharmonic exhaustion function, smooth on $Z_{\mathrm{reg}}$, whose level sets contain at most $p$ points of $D$ each.
\end{lemma}

\begin{proof}
Enumerate the countably many pairs of distinguishable points of $D$ as $(x_k,y_k)$, $k\ge1$, and choose smooth plurisubharmonic functions $u_k$ with $u_k(x_k)\ne u_k(y_k)$. Fix a locally finite atlas of $Z_{\mathrm{reg}}$ and compact sets $K_k$ exhausting $Z_{\mathrm{reg}}$, and choose $\delta_k>0$ such that $\delta_k$ times the maximum of $|u_k|$ on $\{\varphi\le k\}$, and $\delta_k$ times the $C^k$ norm of $u_k$ on $K_k$, are at most $2^{-k}$. For $\varepsilon=(\varepsilon_k)$ in $Q=\prod_k[0,\delta_k]$, the series $v_\varepsilon=\sum_k\varepsilon_ku_k$ converges locally uniformly on $Z$, and with all derivatives uniformly on compact subsets of $Z_{\mathrm{reg}}$. Hence $v_\varepsilon$ is continuous, smooth on $Z_{\mathrm{reg}}$, and plurisubharmonic, since a locally uniform limit of plurisubharmonic functions on a complex space is plurisubharmonic~\cite{FN80}; moreover $|v_\varepsilon|\le C(t)$ on $\{\varphi\le t\}$, where $C(t)=\sum_k\delta_k\max_{\{\varphi\le t\}}|u_k|$ is finite and nondecreasing and does not depend on $\varepsilon$. Choose a smooth convex increasing function $\chi$ with $\chi(t)\ge t+C(t)$. Then $\varphi_\varepsilon=\chi(\varphi)+v_\varepsilon\ge\varphi$ is a continuous plurisubharmonic exhaustion function, smooth on $Z_{\mathrm{reg}}$. Give $Q$ the product of the normalized Lebesgue measures. For each $k$, the set of $\varepsilon\in Q$ with $\varphi_\varepsilon(x_k)=\varphi_\varepsilon(y_k)$ is a null set, since for fixed $(\varepsilon_l)_{l\ne k}$ at most one value of $\varepsilon_k$ satisfies this equation. Hence for almost every $\varepsilon$, distinguishable points of $D$ have distinct values, and the points of $D$ on a level set of $\varphi_\varepsilon$ lie in one class.
\end{proof}

\section{Interpolation by \texorpdfstring{$n$}{n}-forms}\label{sec:interpolation}

\subsection{Proof of Theorem~\ref{thm:main}}
In this subsection, $Y$, $\iota$ and $\mathcal{O}_Y(1)$ are as in Theorem~\ref{thm:main}, and $\mathcal{O}_Y(1)$ carries the pullback of the Fubini--Study metric $h_{\mathrm{FS}}$. Let $C\subset X$ be the set of points at which $d(\zeta\circ\pi)$ has rank less than $n$ for a local embedding $\zeta\colon V\to\mathbb{C}^d$ of an open set $V\subset Y$. This analytic subset does not depend on the local embedding, and over $Y_{\mathrm{reg}}$ it is the set where $\pi$ is not a local biholomorphism, locally the zero set of a Jacobian determinant. It is nowhere dense, since otherwise $\zeta\circ\pi$ would have constant rank less than $n$ on a nonempty open set and the rank theorem would give positive-dimensional fibers. Since $\iota$ is an immersion, $\iota^*\omega_{\mathrm{FS}}$ is a K\"ahler form on $Y_{\mathrm{reg}}$.

\begin{lemma}\label{lem:hyperplane-complement}
Let $A\subset\mathbb{P}^N$ be a hypersurface. Then $Y\setminus\iota^{-1}(A)$ is Stein.
\end{lemma}

\begin{proof}
Let $A$ be the zero set of a homogeneous polynomial $P$ of degree $d$, and let $|\ell|^2=|\ell_0|^2+\dots+|\ell_N|^2$ for unitary homogeneous coordinates. The function $\log(|\ell|^{2d}/|P|^2)$ on $\mathbb{P}^N\setminus A$ is bounded below, tends to infinity at $A$, and has Levi form $d\,\omega_{\mathrm{FS}}>0$. Hence $u=\log(|\ell|^{2d}/|P|^2)\circ\iota$ is a smooth strictly plurisubharmonic function on $Y\setminus\iota^{-1}(A)$ that is bounded below and tends to infinity at $\iota^{-1}(A)$. Let $\varphi$ be a smooth plurisubharmonic exhaustion function of $Y$. Since $u$ is bounded below and tends to infinity at $\iota^{-1}(A)$, the function $\varphi+u$ is a smooth strictly plurisubharmonic exhaustion function of $Y\setminus\iota^{-1}(A)$. Hence $Y\setminus\iota^{-1}(A)$ is Stein~\cite{N62}.
\end{proof}

We prove Theorem~\ref{thm:main} in the following form, which assumes local Steinness only off a hyperplane section and the singular locus; Theorem~\ref{thm:main} is the case of a hyperplane $H$ such that $\iota^{-1}(H)$ contains no irreducible component of $Y$; such hyperplanes exist, since $Y$ has countably many irreducible components.

\begin{theorem}\label{thm:main-refined}
Let $H\subset\mathbb{P}^N$ be a hyperplane such that $\iota^{-1}(H)$ contains no irreducible component of $Y$, and let $\pi\colon X\to Y$ be a nonsingular Riemann domain that is locally Stein over $Y\setminus(\iota^{-1}(H)\cup Y_{\mathrm{sing}})$. Then, with $L$, $k$, $p$, $m$, $S$ and $(a_x)_{x\in S}$ as in Theorem~\textup{\ref{thm:main}}, the conclusion of Theorem~\textup{\ref{thm:main}} holds.
\end{theorem}

\begin{proof}
Put $Y_S=\pi(S)$, a closed discrete subset of $Y$, and $\hat S=\pi^{-1}(Y_S)$, which is closed and discrete. For $y\in Y_S$ let $S_y=S\cap\pi^{-1}(y)$, a finite set, and put $b_y=\iota(y)$ and $\mu=n+k$. Let $h$ be the product of the metric of $L$ and the pullback of $h_{\mathrm{FS}}^m$ on $\mathcal{L}=L\otimes\pi^*\mathcal{O}_Y(m)$; we construct a section of $K_X\otimes\mathcal{L}$.

\emph{The unramified part.} By Lemma~\ref{lem:hyperplane-complement}, $Y_H=Y\setminus\iota^{-1}(H)$ is Stein. Choose $g\in\mathcal{O}(Y_H)$ vanishing on $Y_{H,\mathrm{sing}}$ but not identically on any irreducible component of $Y_H$; such a function exists by Cartan's Theorem~A and the Baire theorem. Put $Y^0=Y_H\setminus g^{-1}(0)$, a Stein manifold, and $X^0=\pi^{-1}(Y^0)\setminus C$. Then $X^0$ is Stein. Indeed, if $U\subset Y^0$ is open with $\pi^{-1}(U)$ Stein, then $\pi^{-1}(U)\setminus C$, locally the complement of the zero set of a holomorphic function in the Stein manifold $\pi^{-1}(U)$, is Stein by Docquier and Grauert~\cite{DG}. Hence each of the countably many connected components of $X^0$ is a locally Stein unramified Riemann domain over the Stein manifold $Y^0$, and is Stein, again by~\cite{DG}. The sets $\iota^{-1}(H)$ and $g^{-1}(0)$ have dimension less than $n$, and a holomorphic map with discrete fibers cannot map a nonempty open subset of the $n$-dimensional manifold $X$ into such a set. Hence $\pi^{-1}(\iota^{-1}(H))$ is a nowhere dense analytic subset of $X$, and $X^0$ is the complement of the nowhere dense analytic subset $(g\circ\pi)^{-1}(0)\cup C$ of $X\setminus\pi^{-1}(\iota^{-1}(H))$. The form $\omega=\frac12\pi^*\iota^*\omega_{\mathrm{FS}}$ is a K\"ahler form on $X^0$.

\emph{Poles.} For $b\in\mathbb{P}^N$ choose unitary homogeneous coordinates $\ell_0,\dots,\ell_N$ with $b=[1:0:\dots:0]$, and put
\[
\psi_{b,\varepsilon}=\log\frac{|\ell_1|^2+\dots+|\ell_N|^2+\varepsilon|\ell_0|^2}{|\ell_0|^2+\dots+|\ell_N|^2}\qquad(0\le\varepsilon\le1).
\]
This function on $\mathbb{P}^N$ satisfies $\psi_{b,\varepsilon}\le0$, it is smooth for $\varepsilon>0$, it decreases to $\psi_{b,0}$ as $\varepsilon\downarrow0$, and $i\partial\bar\partial\psi_{b,\varepsilon}\ge-\omega_{\mathrm{FS}}$, since the logarithm of the numerator is plurisubharmonic on $\mathbb{C}^{N+1}$. In the affine coordinates $\zeta_j=\ell_j/\ell_0$ centered at $b$, we have $\psi_{b,0}=\log\bigl(|\zeta|^2/(1+|\zeta|^2)\bigr)\le\log|\zeta|^2$. On a compact set $Q\subset\mathbb{P}^N\setminus\{b\}$ we have $|\ell_1|^2+\dots+|\ell_N|^2\ge c|\ell|^2$ for some $c>0$; hence $|\psi_{b,\varepsilon}|\le|\psi_{b,0}|$ is bounded on $Q$, and $i\partial\psi_{b,\varepsilon}\wedge\bar\partial\psi_{b,\varepsilon}\le D\,\omega_{\mathrm{FS}}$ on $Q$ with $D$ independent of $\varepsilon\in[0,1]$.

\emph{Localization.} By Lemma~\ref{lem:levels}, $Y$ has a continuous plurisubharmonic exhaustion function $\varphi$, smooth on $Y_{\mathrm{reg}}$, whose level sets contain at most $p$ points of $Y_S$ each. Since $Y_S$ is closed and discrete and the sublevel sets of $\varphi$ are compact, $\varphi(Y_S)$ is a closed discrete subset of $\mathbb{R}$. Let $t_1<t_2<\cdots$ be its elements, finitely or infinitely many, and $G_j=Y_S\cap\varphi^{-1}(t_j)$, a set of at most $p$ points. Choose open intervals $I_j\ni t_j$ with pairwise disjoint closures forming a locally finite family, and $\lambda_j\in C^\infty_c(I_j)$ with $0\le\lambda_j\le1$ and $\lambda_j=1$ near $t_j$. Put
\[
\rho_{j,\varepsilon}=\mu\sum_{y\in G_j}\psi_{b_y,\varepsilon}\circ\iota\le0,\qquad\rho_\varepsilon=\sum_j\lambda_j(\varphi)\,\rho_{j,\varepsilon},
\]
a locally finite sum. The compact set $Z_j=\varphi^{-1}(\operatorname{supp}\lambda_j')$ does not meet $G_j$, and $\iota$ is injective, so there are constants $R_j,D_j$ with $|\rho_{j,\varepsilon}|\le R_j$ and $i\partial\rho_{j,\varepsilon}\wedge\bar\partial\rho_{j,\varepsilon}\le D_j\,\iota^*\omega_{\mathrm{FS}}$ on $Z_j\cap Y_{\mathrm{reg}}$ for $0\le\varepsilon\le1$. On $\varphi^{-1}(I_j)\cap Y_{\mathrm{reg}}$, writing $\lambda_j,\lambda_j',\lambda_j''$ for their values at $\varphi$,
\begin{multline*}
i\partial\bar\partial\rho_\varepsilon=\lambda_j\,i\partial\bar\partial\rho_{j,\varepsilon}+\rho_{j,\varepsilon}\bigl(\lambda_j'\,i\partial\bar\partial\varphi+\lambda_j''\,i\partial\varphi\wedge\bar\partial\varphi\bigr)\\
+\lambda_j'\bigl(i\partial\varphi\wedge\bar\partial\rho_{j,\varepsilon}+i\partial\rho_{j,\varepsilon}\wedge\bar\partial\varphi\bigr).
\end{multline*}
The first term is at least $-p\mu\,\iota^*\omega_{\mathrm{FS}}$. The others vanish off $Z_j$; on $Z_j$ the second is at least $-R_j(|\lambda_j'|\,i\partial\bar\partial\varphi+|\lambda_j''|\,i\partial\varphi\wedge\bar\partial\varphi)$, and by the Cauchy--Schwarz inequality the third is at least $-2D_j\lambda_j'^2\,i\partial\varphi\wedge\bar\partial\varphi-\frac12\iota^*\omega_{\mathrm{FS}}$. Off the sets $\varphi^{-1}(I_j)$, $\rho_\varepsilon=0$.

\emph{Data.} For $y\in Y_S$, choose a relatively compact open neighborhood $P'_y$ of $y$ such that the closures $\overline{P'_y}$ are pairwise disjoint and form a locally finite family, $\overline{P'_y}\cap Y_S=\{y\}$, $\lambda_j=1$ on $\varphi(\overline{P'_y})$ for $y\in G_j$, and $\iota$ maps a neighborhood of $\overline{P'_y}$ isomorphically onto a closed subspace of an open subset of $\mathbb{P}^N$. For $x\in S$ with $\pi(x)=y$, apply Lemma~\ref{lem:finrep} with $P'=P'_y$ and $W'$ a coordinate ball centered at $x$, chosen pairwise disjoint for the finitely many points of $S_y$. Shrinking $W'$ as in the last part of the lemma, we obtain relatively compact open sets $W_x\ni x$ and one open neighborhood $P_y\subset P'_y$ of $y$ such that $\pi|_{W_x}\colon W_x\to P_y$ is proper and $\overline{W_x}\cap\pi^{-1}(y)=\{x\}$ for $x\in S_y$. Then $\overline{W_x}\cap\hat S=\{x\}$, and the sets $\overline{W_x}$, $x\in S$, are pairwise disjoint and form a locally finite family. Let $A_x$ be a holomorphic section of $K_X\otimes\mathcal{L}$ on a neighborhood of $\overline{W_x}$ with $j^k_xA_x=a_x$, and let $\theta_y$ be a smooth function on $Y$ with compact support in $P_y$, equal to $\Theta_y\circ\iota$ near $\overline{P_y}$ for some $\Theta_y\in C^\infty(\mathbb{P}^N)$, such that $0\le\theta_y\le1$ and $\theta_y=1$ on a neighborhood of $y$. Define
\[
\tilde s=\sum_{x\in S}(\theta_{\pi(x)}\circ\pi)\,A_x,
\]
the $x$th term being taken on $W_x$ and extended by $0$. By the remark after Lemma~\ref{lem:finrep}, each term is smooth with compact support, and the sum is locally finite. Thus $\tilde s$ is a smooth section of $K_X\otimes\mathcal{L}$. It equals $A_x$ near $x\in S$, and it vanishes near $\hat S\setminus S$; in particular it is holomorphic near $\hat S$. Its derivative
\[
v=\bar\partial\tilde s=\sum_{x\in S}\pi^*(\bar\partial\theta_{\pi(x)})\wedge A_x
\]
is supported in the union $T$ of the compact sets $T_x=(\pi|_{W_x})^{-1}(\operatorname{supp} d\theta_{\pi(x)})$, which are disjoint from $\hat S$. On $Y_{\mathrm{reg}}$, a covector restricted to a subspace does not increase its norm, so $|\bar\partial\theta_y|_{\iota^*\omega_{\mathrm{FS}}}\le\sup|\bar\partial\Theta_y|_{\omega_{\mathrm{FS}}}$. Since $\pi$ is a local isometry from $(X^0,\omega)$ to $(Y_{\mathrm{reg}},\frac12\iota^*\omega_{\mathrm{FS}})$, \eqref{eq:product} gives, for $x\in S_y$,
\begin{equation}\label{eq:cutoffnorm}
|v|^2_{\omega,h}\,dV_\omega=2\bigl(|\bar\partial\theta_y|^2_{\iota^*\omega_{\mathrm{FS}}}\circ\pi\bigr)\,|A_x|^2_h\le2\sup|\bar\partial\Theta_y|^2_{\omega_{\mathrm{FS}}}\,|A_x|^2_h\quad\text{on }W_x\cap X^0,
\end{equation}
a bound that does not involve $\omega$, which degenerates outside $X^0$.

\emph{Weights.} Let $x\in S_y$ and $y\in G_j$. On $T_x$ we have $\rho_\varepsilon\circ\pi=\rho_{j,\varepsilon}\circ\pi\ge\rho_{j,0}\circ\pi\ge-c_y$ for a constant $c_y$, since $\iota\circ\pi$ maps $T_x$ into the complement of the points $b_{y'}$, $y'\in G_j$. Put $t_y=\min_{\operatorname{supp}\theta_y}\varphi$ and
\[
M_y=2e^{c_y}\sum_{x\in S_y}\sup|\bar\partial\Theta_y|^2_{\omega_{\mathrm{FS}}}\int_{T_x}|A_x|^2_h ,
\]
and enumerate $Y_S$ by an injective map $\nu\colon Y_S\to\{1,2,\dots\}$. Since only finitely many $t_y$ lie below any given number, there is a smooth convex increasing function $\kappa$ on $\mathbb{R}$ such that
\begin{equation}\label{eq:kappa}
\begin{aligned}
\kappa'&\ge R_j|\lambda_j'|,\quad\kappa''\ge R_j|\lambda_j''|+2D_j\lambda_j'^2\quad\text{on }I_j\qquad(j\ge1),\\
\kappa(t_y)&\ge\log\bigl(1+2^{\nu(y)}M_y\bigr)\qquad(y\in Y_S).
\end{aligned}
\end{equation}
For $\varepsilon\in[0,1]$ put $\Psi_\varepsilon=(\kappa(\varphi)+\rho_\varepsilon)\circ\pi$, which is smooth on $X^0$ for $\varepsilon>0$. Since $i\partial\bar\partial\kappa(\varphi)=\kappa'(\varphi)\,i\partial\bar\partial\varphi+\kappa''(\varphi)\,i\partial\varphi\wedge\bar\partial\varphi$, the estimates above give $i\partial\bar\partial(\kappa(\varphi)+\rho_\varepsilon)\ge-(p\mu+\frac12)\,\iota^*\omega_{\mathrm{FS}}$ on $Y_{\mathrm{reg}}$, hence
\[
i\Theta(h)+i\partial\bar\partial\Psi_\varepsilon\ge\bigl(m-p\mu-\tfrac12\bigr)\,\pi^*\iota^*\omega_{\mathrm{FS}}\ge\omega\quad\text{on }X^0 .
\]
On $T_x$ we have $\kappa(\varphi\circ\pi)\ge\kappa(t_y)$, so~\eqref{eq:cutoffnorm} and~\eqref{eq:kappa} give
\[
\int_{X^0}|v|^2_{\omega,h}\,dV_\omega\,e^{-\Psi_\varepsilon}\le\sum_{y\in Y_S}e^{-\kappa(t_y)}M_y\le\sum_{y\in Y_S}2^{-\nu(y)}\le1\qquad(0<\varepsilon\le1).
\]

\emph{Solving.} By Lemma~\ref{lem:L2} on the Stein manifold $X^0$, there are $\mathcal{L}$-valued $(n,0)$-forms $u_\varepsilon$ on $X^0$ with $\bar\partial u_\varepsilon=v$ and $\int_{X^0}|u_\varepsilon|^2_he^{-\Psi_\varepsilon}\le1$. Since $\Psi_\varepsilon\le\Psi_{\varepsilon'}$ for $\varepsilon\le\varepsilon'$, a diagonal subsequence of $(u_{2^{-l}})$ converges weakly in the weighted $L^2$ space defined by $h$ and $e^{-\Psi_{2^{-l'}}}$ for every integer $l'\ge1$. The limits agree as distributions, since these weights are smooth and positive on $X^0$. The common limit $u$ satisfies $\bar\partial u=v$, and weak lower semicontinuity followed by monotone convergence gives
\begin{equation}\label{eq:weightedn}
\int_{X^0}|u|^2_{h}e^{-\Psi_0}\le1 .
\end{equation}
The function $\Psi_0\le\kappa(\varphi)\circ\pi$ is bounded above on compact subsets of $X$. As $h$ is smooth on $X$, \eqref{eq:weightedn} shows that the coefficient of $u$ in a local frame of $K_X\otimes\mathcal{L}$ is square integrable near every point of $X$. Hence the section $\tilde s-u$, holomorphic on $X^0$, extends across $(g\circ\pi)^{-1}(0)\cup C$ to $X\setminus\pi^{-1}(\iota^{-1}(H))$, and then across $\pi^{-1}(\iota^{-1}(H))$, to a section $s\in H^0(X,K_X\otimes\mathcal{L})$~\cite[Chapter~VIII, Lemma~7.3]{Dem}.

\emph{Jets.} Let $z\in\hat S$, $y=\pi(z)\in G_j$, and choose coordinates $w$ centered at $z$ on a ball $U$ on which $\tilde s$ is holomorphic. Near $z$ we have $\lambda_j(\varphi\circ\pi)=1$, the other functions $\lambda_{j'}(\varphi\circ\pi)$ vanish, and the terms of $\rho_{j,0}$ are nonpositive, so $\Psi_0\le c+\mu\,\psi_{b_y,0}\circ\iota\circ\pi\le c'+\mu\log|w|^2$, since the affine coordinates centered at $b_y$, composed with $\iota\circ\pi$, vanish at $z$. By~\eqref{eq:weightedn}, the coefficient of $u=\tilde s-s$ in a frame over $U$ is square integrable against $|w|^{-2\mu}d\lambda$ near $0$, and Lemma~\ref{lem:jets} gives $u_z\in\mathfrak{m}_z^{\mu-n+1}(K_X\otimes\mathcal{L})_z$, that is, $j^k_zu=0$. For $x\in S_y$ we have $\tilde s=A_x$ near $x$, so $j^k_xs=j^k_xA_x=a_x$. For $z\in\pi^{-1}(y)\setminus S$ we have $\tilde s=0$ near $z$, so $j^k_zs=0$.
\end{proof}

\subsection{Consequences}\label{sec:consequences}

\begin{proof}[Proof of the Stein case of Theorem~\textup{\ref{thm:main}}]
Let $F\colon Y\to\mathbb{C}^{2n+1}$ be proper and injective~\cite{N60}. Then $Y'=F(Y)$, with its reduced structure, is a closed analytic subset of $\mathbb{C}^{2n+1}$, hence a Stein space, and $F$ is a homeomorphism onto $Y'$; it need not be an immersion, since Theorem~\ref{thm:main} will be applied over $Y'$. The map $\pi'=F\circ\pi$ is a nonsingular Riemann domain over $Y'$ such that $S$ is finite over $Y'$ and $\pi'^{-1}(\pi'(S))=\pi^{-1}(\pi(S))$. It is locally Stein: if $U\ni y$ is open with $\pi^{-1}(U)$ Stein and $B$ is a ball centered at $F(y)$ with $B\cap Y'\subset F(U)$, then $\pi'^{-1}(B)$ is the preimage of $B$ under $\pi'|_{\pi^{-1}(U)}$, hence Stein. The classes of indistinguishable points of $Y'$ are points, and for the inclusion $\iota$ of $Y'$ in $\mathbb{C}^{2n+1}\subset\mathbb{P}^{2n+1}$, division by a power of $\ell_0\circ\iota\circ\pi'$, where $\ell_0$ defines the hyperplane at infinity, identifies $K_X\otimes L\otimes\pi'^*\mathcal{O}_{Y'}(m)$ with $K_X\otimes L$. Now apply Theorem~\ref{thm:main} with $p=1$ and $m=n+k+1$.
\end{proof}

\begin{lemma}\label{lem:genpos}
Let $\pi\colon X\to Y$ be a holomorphic map with discrete fibers from a complex manifold of dimension $n$ to a complex space, and let $L$ be a holomorphic line bundle on $X$ such that on $V=H^0(X,L)$ the maps $s\mapsto j^1_xs$ and $s\mapsto(j^0_{x_1}s,j^0_{x_2}s)$ are onto for all $x\in X$ and all $x_1\ne x_2$ with $\pi(x_1)=\pi(x_2)$. Then the tuples $s\in V^{n+2}$ without common zero for which $(\pi,[s_0:\dots:s_{n+1}])\colon X\to Y\times\mathbb{P}^{n+1}$ is an injective immersion form a dense $G_\delta$ subset of $V^{n+2}$.
\end{lemma}

\begin{proof}
Put $A=(X\times_YX)\setminus\Delta_X$, an analytic subset of $(X\times X)\setminus\Delta_X$. For $x\in X$ let $d\pi_x$ be the differential of $\zeta\circ\pi$ at $x$ for a local embedding $\zeta$ of $Y$ near $\pi(x)$; its kernel does not depend on $\zeta$, since the components of another local embedding are holomorphic functions of $\zeta$ near $\pi(x)$. For $1\le r\le n$ let $C_r=\{x\in X:\operatorname{rank}d\pi_x=n-r\}$, a locally closed analytic subset, and $C=\bigcup_{r\ge1}C_r$. We claim that $\dim A\le n$ and $\dim C_r\le n-r$. Let $\Sigma$ be a connected complex submanifold of dimension $d$ contained in $A$, respectively in $C_r$, and let $\psi$ be the restriction to $\Sigma$ of the first projection, respectively of $\pi$. Then $\psi$ has discrete fibers, so by the rank theorem its rank is $d$ on the dense open set where it is maximal; on the other hand, this rank is at most $n$, respectively $n-r$. Writing $A$ and $C_r$ as unions of the regular loci of their iterated singular loci, we obtain countably many locally closed submanifolds $A_\alpha$ and $\Sigma_\beta$ with $\dim A_\alpha\le n$ and $\dim\Sigma_\beta\le n-r$ covering $A$ and $C_r$. Over each $\Sigma_\beta$, $\ker d\pi$ is a holomorphic subbundle of rank $r$ of $T_X|_{\Sigma_\beta}$, being the kernel of a bundle map of constant rank.

In local frames of $L$, a tuple $s$ is a map $\sigma$ into $\mathbb{C}^{n+2}$. The map $(\pi,[\sigma])$ is an injective immersion if and only if: (1) $\sigma(x)\ne0$ for all $x$; (2) for $(x_1,x_2)\in A$, the vectors $\sigma(x_1),\sigma(x_2)$ are linearly independent; (3) for $x\in C_r$ and a basis $\xi_1,\dots,\xi_r$ of $\ker d\pi_x$, the vectors $\sigma(x),\xi_1\sigma(x),\dots,\xi_r\sigma(x)$ are linearly independent. Indeed, points with different images are separated by $\pi$, off $C$ the map $\pi$ is immersive, and $d[\sigma]_x(\xi)=0$ exactly when $\xi\sigma(x)\in\mathbb{C}\sigma(x)$. A change of frame multiplies $\sigma$ by a nonvanishing function $\lambda$ and adds $(\xi\lambda)\sigma$ to $\xi\sigma$, so these conditions do not depend on the frames. The bad sets in the respective jet spaces are cones of codimension $n+2$ for (1), $n+1$ for (2), and $n-r+2$ for (3); each exceeds the dimension of the corresponding space $X$, $A_\alpha$ or $\Sigma_\beta$.

Let $Q$ be one of $X$, $A_\alpha$, $\Sigma_\beta$, with the corresponding condition, and let $K\subset Q$ be compact. The set $G(K)\subset V^{n+2}$ of tuples satisfying the condition at every point of $K$ is open, since locally uniform convergence of holomorphic sections implies $C^1$ convergence on compact sets. To prove density, fix $s\in V^{n+2}$. The hypotheses on $V$ and compactness give finitely many $q_1,\dots,q_P\in V$ that span the required values or jets near $K$: values for (1), values at both points for (2), and values together with derivatives along $\ker d\pi$ for (3). For $t=(t_{l,i})\in\mathbb{C}^{P(n+2)}$ put $s^t_i=s_i+\sum_lt_{l,i}q_l$. In local frames, the map from $(x,t)$ to these values or jets is a submersion, since it is affine and surjective in $t$. Stratify each determinantal cone by rank. Since its codimension exceeds $\dim Q$, each inverse image stratum has complex dimension less than $P(n+2)$, so its projection to $\mathbb{C}^{P(n+2)}$ has measure zero. A countable coordinate cover gives the same conclusion for all parameters that fail the condition somewhere near $K$, so $s^t\in G(K)$ for arbitrarily small $t$.

$V^{n+2}$ is a Fr\'echet space, and $X$, $A_\alpha$ and $\Sigma_\beta$ are countable unions of compact sets. The intersection of the countably many open dense sets $G(K)$ is a dense $G_\delta$, and its elements satisfy (1)--(3) everywhere.
\end{proof}

\begin{proof}[Proof of Corollary~\textup{\ref{cor:main}}]
(i) Take $L$ trivial in Theorem~\ref{thm:main}. The case $k=0$, $p=1$ shows that $K_X\otimes\pi^*\mathcal{O}_Y(n+1)$ is generated by global sections. For $m\ge n+2$, the case $k=1$, $p=1$ and the case $k=0$, $p=1$ for two points with the same image show that $V=H^0(X,K_X\otimes\pi^*\mathcal{O}_Y(m))$ satisfies the hypotheses of Lemma~\ref{lem:genpos}. If $Y$ is a Stein space, the Stein case of Theorem~\ref{thm:main} with $L$ trivial shows the same for $V=H^0(X,K_X)$. (ii) If $K_X^{-1}$ carries a smooth hermitian metric of semipositive curvature, apply the Stein case of Theorem~\ref{thm:main} with $L=K_X^{-1}$; then $K_X\otimes L$ is trivial. If $\pi$ is unramified, then $\pi(X)\subset Y_{\mathrm{reg}}$, since $X$ is nonsingular, and $K_X=\pi^*K_{Y_{\mathrm{reg}}}$. By Theorem~A, countably many holomorphic vector fields on $Y$ span the tangent space at every point of $Y_{\mathrm{reg}}$, so their $n$-fold wedge products $s_I$ have no common zero on $Y_{\mathrm{reg}}$. Writing $s_I=f_Ie$ in local frames $e$ and choosing $\varepsilon_I>0$ so small that $\sum_I\varepsilon_I|f_I|^2$ converges locally uniformly, hence smoothly, we obtain a smooth metric $|e|^2=(\sum_I\varepsilon_I|f_I|^2)^{-1}$ on $K_{Y_{\mathrm{reg}}}^{-1}$ with semipositive curvature $i\partial\bar\partial\log\sum_I\varepsilon_I|f_I|^2$; its pullback to $K_X^{-1}$ reduces this case to the first.
\end{proof}

\begin{proof}[Proof of Corollary~\textup{\ref{cor:envelope}}]
If $X$ is Stein, then $(X,\pi,\mathrm{id})$ is an envelope of holomorphy. Conversely, let $(\tilde X,p,\tau)$ be one. By Corollary~\ref{cor:main}(ii), holomorphic functions separate the points of $X$. Since $p\circ\tau=\pi$, the map $\tau$ is locally biholomorphic, and since $\tilde f\circ\tau=f$, it is injective, so it maps $X$ biholomorphically onto an open subset $U$ of $\tilde X_{\mathrm{reg}}$. Every holomorphic function on $U$ extends to $\tilde X$.

The set $U$ is locally Stein in $\tilde X$. Indeed, let $q\in\tilde X$, choose an open neighborhood $V$ of $p(q)$ with $\pi^{-1}(V)$ Stein, and choose a Stein open neighborhood $W\subset p^{-1}(V)$ of $q$. Then $U\cap W=\tau(\pi^{-1}(V))\cap W$ is the intersection of two Stein open subsets of $\tilde X$, hence Stein.

The space $\tilde X$ is irreducible. Since $U$ is connected and lies in $\tilde X_{\mathrm{reg}}$, it lies in a single irreducible component. If there were another component, Theorem~B would give a nonzero holomorphic function on $\tilde X$ vanishing on all the other components, hence on $U$. This would contradict the uniqueness of extensions.

Let $\nu\colon\hat X\to\tilde X$ be the normalization. Then $\hat X$ is a connected normal Stein space, and $\nu$ maps $\hat U=\nu^{-1}(U)$ biholomorphically onto $U$. The set $\hat U$ is locally Stein in $\hat X$, since $\hat U\cap\nu^{-1}(W)$ is biholomorphic to $U\cap W$. Every holomorphic function on $\hat U$ extends to $\hat X$.

Suppose a sequence in $\hat U$ converges to a point of $\partial\hat U\setminus\hat X_{\mathrm{sing}}$. Then Preda's theorem~\cite{Pre15} gives a holomorphic function on $\hat U$ that is unbounded on this sequence. This is impossible, since the function extends to $\hat X$. Hence $\partial\hat U\subset\hat X_{\mathrm{sing}}$. Since $\hat X_{\mathrm{reg}}$ is connected~\cite[Chapter~7, Section~4.2]{GR84}, it follows that $\hat U=\hat X_{\mathrm{reg}}$.

Finally, let $s\in\hat X_{\mathrm{sing}}$, and let $W'$ be an open neighborhood of $s$ such that $\hat U\cap W'=W'\setminus\hat X_{\mathrm{sing}}$ is Stein. By the Riemann extension theorem on normal spaces~\cite[Chapter~7, Section~4.2]{GR84}, every holomorphic function on $W'\setminus\hat X_{\mathrm{sing}}$ extends to $W'$. So every such function is bounded on a sequence converging to $s$. This contradicts holomorphic convexity. Hence $\hat X=\hat U$ is nonsingular, and $X\cong\hat U=\hat X$ is Stein.
\end{proof}

\subsection{Interpolation and Steinness}\label{sec:steinness}
In this subsection, $L$ is a holomorphic line bundle with a smooth hermitian metric $h$ on a complex manifold $X$, and its sections interpolate arbitrary values on every closed discrete subset of $X$. The proof of Proposition~\ref{prop:steinness} passes through Grauert's convexity with respect to line bundles~\cite[p.~380]{Gr63}, which was studied by V\^aj\^aitu~\cite{Vaj16} and by the second author~\cite{Oh24,Oh24b}.

\begin{lemma}\label{lem:Kstar}
Every compact set $K\subset X$ has a compact neighborhood $K^\star$ such that for every $x\in X\setminus K^\star$ and every $c>0$ there is $s\in H^0(X,L)$ with $|s(x)|>c\sup_K|s|$.
\end{lemma}

\begin{proof}
Suppose not. Taking compact neighborhoods of $K$ that exhaust $X$, we find distinct points $x_\nu$, $\nu\ge1$, leaving every compact subset of $X$, and constants $c_\nu>0$ such that $|s(x_\nu)|\le c_\nu\sup_K|s|$ for all $s\in H^0(X,L)$. The points $x_\nu$ form a closed discrete set. Choose $w_\nu$ in the fiber of $L$ at $x_\nu$ with $|w_\nu|=\nu c_\nu$. A section $s$ with $s(x_\nu)=w_\nu$ for all $\nu$ would satisfy $\sup_K|s|\ge\nu$ for every $\nu$, which is impossible.
\end{proof}

For each compact set $K\subset X$ we fix such a neighborhood $K^\star$.

\begin{lemma}\label{lem:stein-neighborhoods}
Let $X$ carry a smooth strictly plurisubharmonic function $\psi$. Every compact subset of $X$ is contained in a relatively compact Stein open subset of $X$.
\end{lemma}

\begin{proof}
Let $K\subset X$ be a nonempty compact set, and choose a relatively compact open set $U$ containing $K^\star$, and a constant $A>0$ such that $A\,i\partial\bar\partial\psi-i\Theta(h)>0$ on a neighborhood of $\overline U$. For every nonzero $s\in H^0(X,L)$, the function $A\psi+\log|s|_h^2$ is strictly plurisubharmonic there. Put $b=A\max_K\psi+1$, and choose $M>0$ with $A\min_{\partial U}\psi+M>b$. Since $\partial U\cap K^\star=\emptyset$, Lemma~\ref{lem:Kstar} gives, for every $x\in\partial U$, some $s\in H^0(X,L)$ with $|s|_h<1$ on $K$ and $\log|s(x)|_h^2>M$, after rescaling. By compactness of $\partial U$, finitely many such sections $s_1,\dots,s_r$ give a continuous strictly plurisubharmonic function
\[
v=A\psi+\max\bigl(0,\log|s_1|_h^2,\dots,\log|s_r|_h^2\bigr)
\]
on a neighborhood of $\overline U$, with $v<b$ on $K$ and $v>b$ on $\partial U$. The open set $D=\{x\in U:v(x)<b\}$ contains $K$ and is relatively compact in $U$. The function $(b-v)^{-1}$ is a continuous strictly plurisubharmonic exhaustion of $D$, so $D$ is Stein~\cite{Gr58,N62}.
\end{proof}

\begin{lemma}\label{lem:stein-hull-comparison}
Let $D\subset X$ be a Stein open subset, and let $K\subset D$ be compact. Then $\widehat K_{\mathcal{O}(D)}\subset K^\star$.
\end{lemma}

\begin{proof}
Let $x\in\widehat K_{\mathcal{O}(D)}$. Cartan's Theorem~A on $D$ gives $t\in H^0(D,L^{-1})$ with $t(x)\ne0$. For every $s\in H^0(X,L)$, the product $(s|_D)t$ is holomorphic on $D$, so $|s(x)||t(x)|\le\sup_K|s|\sup_K|t|$. By Lemma~\ref{lem:Kstar}, $x\in K^\star$.
\end{proof}

\begin{proof}[Proof of Proposition~\textup{\ref{prop:steinness}}]
If $X$ is Stein, the sections of every holomorphic line bundle interpolate arbitrary values on closed discrete sets, by Cartan's Theorem~B. Conversely, let $L$ be as above. Since the Stein space $Y$ carries a smooth strictly plurisubharmonic function, so does $X$, by the case $q=1$ of a theorem of V\^aj\^aitu~\cite[Theorem~C]{Vaj94}. Given a compact set $A\subset X$, Lemma~\ref{lem:stein-neighborhoods} gives compact sets $K_j$ exhausting $X$ and relatively compact Stein open sets $D_j$ such that
\[
A\subset K_1,\qquad K_j^\star\subset D_j\Subset\mathring K_{j+1}\quad(j\ge1):
\]
start with a compact $K_1\supset A$ and apply the lemma to $K_1^\star$; inductively, choose $K_{j+1}$ whose interior contains $\overline D_j$ and the next set of a fixed compact exhaustion, and apply the lemma to $K_{j+1}^\star$. By Lemma~\ref{lem:stein-hull-comparison}, $\widehat{K_j}_{\mathcal{O}(D_{j+1})}\subset K_j^\star\subset D_j$. Consequently, every $f_j\in\mathcal{O}(D_j)$ can be approximated uniformly on $K_j$ by functions in $\mathcal{O}(D_{j+1})$, by Oka--Weil in the Stein manifold $D_{j+1}$. For $f_1\in\mathcal{O}(D_1)$ and $\varepsilon>0$, choose successively $f_{j+1}\in\mathcal{O}(D_{j+1})$ with $\sup_{K_j}|f_{j+1}-f_j|<2^{-j}\varepsilon$. The functions $f_j$ converge uniformly on compact subsets of $X$ to $f\in\mathcal{O}(X)$, with $\sup_A|f-f_1|\le\varepsilon$. Thus $\mathcal{O}(X)|_A$ is dense in $\mathcal{O}(D_1)|_A$.

Taking $A$ to contain two prescribed points shows that global holomorphic functions separate them, since $D_1$ is Stein, and taking $A$ to contain a closed coordinate ball shows that they give local coordinates, by Cauchy estimates. Finally, fix a compact set $K\subset X$ and $x\notin K^\star$, and apply the construction with $A=K\cup\{x\}$. By Lemma~\ref{lem:stein-hull-comparison}, $x\notin\widehat K_{\mathcal{O}(D_1)}$, so some $f_1\in\mathcal{O}(D_1)$ satisfies $|f_1(x)|>\sup_K|f_1|$, and approximation on $A$ gives a function in $\mathcal{O}(X)$ with the same property. Hence $\widehat K_{\mathcal{O}(X)}\subset K^\star$ is compact, $X$ is holomorphically convex, and therefore Stein.
\end{proof}

\section{Forn\ae ss' surface}\label{sec:F78}

\subsection{Holomorphic functions}
The following elementary lemma isolates the continuation argument of~\cite[Section~2, Lemma~7]{F77}, which was also used by Vigna Suria~\cite{VS84}. We then verify its hypotheses for Forn\ae ss' ramified surface of 1978.

\begin{lemma}\label{lem:descent}
Let $X$ be a connected complex manifold and $\pi\colon X\to\Delta\times\mathbb{C}$ a holomorphic map. Suppose that, for some $R_0>0$, $\pi$ maps $\pi^{-1}(\Delta\times\{|w|>R_0\})$ biholomorphically onto $\Delta\times\{|w|>R_0\}$. Suppose also that there are distinct points $a_n\in\Delta$ converging to some $a\in\Delta$, and holomorphic maps $\sigma_n\colon\mathbb{C}\to X$ with $\pi\circ\sigma_n(w)=(a_n,w)$. Then $\pi^*\colon\mathcal{O}(\Delta\times\mathbb{C})\to\mathcal{O}(X)$ is an isomorphism.
\end{lemma}

\begin{proof}
Given $f\in\mathcal{O}(X)$, use the exterior biholomorphism to write $f$ as a function $g(z,w)$ for $|w|>R_0$. For $R>R_0$, the Laurent coefficients
\[
c_j(z)=\frac1{2\pi i}\int_{|w|=R}g(z,w)w^{-j-1}\,dw
\]
are holomorphic on $\Delta$ and independent of $R$. Since $g(a_n,w)=f(\sigma_n(w))$ for $|w|>R_0$ and $f\circ\sigma_n$ is entire, we have $c_j(a_n)=0$ for every $j<0$. By the identity theorem, $c_j=0$ for $j<0$. Cauchy's estimates, applied with arbitrarily large $R$, show that $h(z,w)=\sum_{j\ge0}c_j(z)w^j$ converges locally uniformly on $\Delta\times\mathbb{C}$. Thus $f=h\circ\pi$ on the exterior chart and, by connectedness, on all of $X$. The injectivity of $\pi^*$ also follows from the exterior chart.
\end{proof}

We recall the construction of~\cite{F78}. Let $A_1=\{|w|>1\}$. By~\cite[Lemma~1]{F78}, proved in~\cite[Section~2]{F77}, there are a negative subharmonic function $\vartheta$ on $\Delta$, pairwise disjoint disks $\Delta_n\subset\Delta$ centered at $a_n=1/n$ $(n\ge2)$, smaller concentric disks $\Delta_n'\subset\Delta_n$, and integers $m_n\ge1$ such that
\[
\vartheta>-1/2\quad\text{off }\bigcup\Delta_n,\qquad
\vartheta=m_n^{-1}\log|(z-a_n)/2|-1\quad\text{on }\Delta_n'.
\]
Put
\[
\begin{aligned}
\Omega_\vartheta&=\{\vartheta(z)<\log|w|\}\subset\Delta\times\mathbb{C}^*,\qquad
U_n=\{|\eta|<2e^{m_n}\}\subset\mathbb{C}^2,\\
\Psi_n(\eta,w)&=(a_n+C_n\eta w^{1+m_n}+\varepsilon_n\eta^2,w).
\end{aligned}
\]
Let $\Delta_n''$ be the disk concentric with $\Delta_n'$ of half its radius. As in~\cite[Lemma~2]{F78}, first choose $C_n\ge1$ large, then disks $\Delta_n^4\Subset\Delta_n^3\Subset\Delta_n''$ containing $a_n$, and finally $\varepsilon_n>0$ small. Then $\Psi_n$ is injective and locally biholomorphic on $\Gamma_n=U_n\cap\Psi_n^{-1}((\Delta_n'\setminus\overline{\Delta_n^4})\times\mathbb{C})$ by~\cite[Lemma~2(i)--(ii)]{F78}. The surface $(M,\pi)$ is glued from three kinds of pieces:
\begin{itemize}
\item over $\Delta\setminus\bigcup\overline{\Delta_n''}$, the domain $\Omega_\vartheta$ with the inclusion;
\item over $\Delta_n'\setminus\overline{\Delta_n^4}$, the set $\Omega_\vartheta\cap\Psi_n(\Gamma_n)$ with the inclusion;
\item over $\Delta_n^3$, the set $V_n=U_n\cap\Psi_n^{-1}(\Delta_n^3\times\mathbb{C})$ with $\pi=\Psi_n$.
\end{itemize}
The first two descriptions agree over $\Delta_n'\setminus\overline{\Delta_n''}$ by~\cite[Lemma~2(iii)]{F78}, and the last two over $\Delta_n^3\setminus\overline{\Delta_n^4}$ by~\cite[Lemma~2(iv)]{F78}; the pieces are glued accordingly. Viewed as a map to $\mathbb{C}^2$, $\pi$ is locally Stein by~\cite[Lemma~4]{F78}, and its image lies in $\Delta\times\mathbb{C}$. Forn\ae ss also proves that $M$ is not Stein; this also follows from Theorem~\ref{thm:fornaess}(a). In accordance with our definition of a Riemann domain, we regard $M$ as connected, replacing it if necessary by its connected component containing $\pi^{-1}(\Delta\times A_1)$, which is again locally Stein over $\Delta\times\mathbb{C}$. Since $\varepsilon_n$ is chosen sufficiently small after $C_n$ and the disks, we have, without any further restriction on the construction,
\begin{equation}\label{eq:eps}
4\varepsilon_ne^{m_n}<C_n,\qquad\{|z-a_n|\le4\varepsilon_ne^{2m_n}\}\subset\Delta_n^3.
\end{equation}

\begin{lemma}\label{lem:fornaess-local}
\begin{enumerate}[label=\textup{(\alph*)}]
\item The map $\pi\colon\pi^{-1}(\Delta\times A_1)\to\Delta\times A_1$ is biholomorphic.
\item The map $\sigma_n(w)=(0,w)\in V_n$ satisfies $\pi\circ\sigma_n(w)=(a_n,w)$.
\item For $0<|\eta|<2e^{m_n}$, the points $(\pm\eta,0)\in V_n$ are distinct and have the same image.
\item The points in~\textup{(b)} and~\textup{(c)} lie in the connected component of $M$ containing $\pi^{-1}(\Delta\times A_1)$.
\end{enumerate}
\end{lemma}

\begin{proof}
(a) Since $\vartheta<0$, we have $\Delta\times A_1\subset\Omega_\vartheta$. Over $\Delta_n'$, the lifts of $(z,w)$ correspond to the roots $\eta$ with $|\eta|<2e^{m_n}$ of $\varepsilon_n\eta^2+C_nw^{1+m_n}\eta=z-a_n$. For $|\eta|=2e^{m_n}$ and $|w|>1$, \eqref{eq:eps} and $|z-a_n|<1<C_ne^{m_n}$ give
\[
|\varepsilon_n\eta^2-(z-a_n)|\le4\varepsilon_ne^{2m_n}+|z-a_n|<2C_ne^{m_n}<|C_nw^{1+m_n}\eta|,
\]
so by Rouch\'e's theorem there is exactly one root, counted with multiplicity. This root is simple and depends holomorphically on $(z,w)$; the gluing identifies it with the inclusion chart. (b) and (c) follow immediately from~\eqref{eq:eps}.

(d) The set $M_\infty=\pi^{-1}(\Delta\times A_1)$ is connected by~(a). The curve $\sigma_n(\mathbb{C})$ meets $M_\infty$ and contains $(0,0)$, and the disk $\{(t,0):|t|<2e^{m_n}\}$ lies in $V_n$ by~\eqref{eq:eps} and joins $(0,0)$ to $(\pm\eta,0)$.
\end{proof}

\begin{proof}[Proof of Theorem~\textup{\ref{thm:fornaess}(a)}]
By Lemma~\ref{lem:fornaess-local}(a), (b) and~(d), Lemma~\ref{lem:descent} applies with $R_0=1$, $a_n=1/n$ and the maps $\sigma_n$ of Lemma~\ref{lem:fornaess-local}(b). By Lemma~\ref{lem:fornaess-local}(c) and~(d), $\pi$ is not injective, so $\mathcal{O}(M)=\pi^*\mathcal{O}(\Delta\times\mathbb{C})$ does not separate points.
\end{proof}

By Corollary~\ref{cor:main}(ii) and Theorem~\ref{thm:fornaess}(a), $K_M^{-1}$ carries no smooth semipositive metric.

\subsection{Increasing unions}
The anticanonical bundle yields a holomorphically separable increasing union. For $k\ge2$, let $\varrho_k<-1$ denote the value of the function $m_k^{-1}\log|(z-a_k)/2|-1$ on $\partial\Delta'_k$; by upper semicontinuity, $\vartheta\ge\varrho_k$ on $\partial\Delta'_k$. Put
\[
M^{(n)}=M\setminus\pi^{-1}\Bigl(\bigcup_{k>n}\{(z,w):z\in\Delta'_k,\ |w|\le\exp(\varrho_k)\}\Bigr).
\]
These sets are open: the closure in $\Delta\times\mathbb{C}$ of the removed set adds only points over $\partial\Delta'_k$ with $|w|\le\exp(\varrho_k)$ and points over $z=0$ with $|w|\le e^{-1}$, and these do not lie in $\pi(M)$, since over them $M$ is the domain $\Omega_\vartheta$, while $\vartheta\ge\varrho_k$ on $\partial\Delta'_k$ and $\vartheta(0)>-1/2$ because $0$ lies in no $\Delta_n$. The sets $M^{(n)}$ increase, and their union is $M$. Forn\ae ss states that $M$ can be exhausted by Stein open subsets~\cite[p.~275]{F78}: by the argument of~\cite[Lemma~4]{F78}, each $M^{(n)}$ carries a strictly plurisubharmonic exhaustion function~\cite[p.~277]{F78}. We use this assertion.

\begin{lemma}[{\cite[Theorems~4 and~5]{MM60}}]\label{lem:punctured-line}
Let $L$ be a holomorphic line bundle on a complex manifold $B$, and let $L^\times$ be its total space with the zero section removed. Then $L^\times$ is Stein if and only if $B$ is Stein.
\end{lemma}

Indeed, $L^\times$ is a holomorphic principal bundle over $B$ with structure group $\mathbb{C}^*$, which is Stein and is the complexification of its maximal compact subgroup.

\begin{proof}[Proof of Theorem~\textup{\ref{thm:fornaess}(b)}]
Let $X$ be the total space of $K_M^{-1}$ with the zero section removed, and $q\colon X\to M$ the projection. The sets $q^{-1}(M^{(n)})$ are Stein by Lemma~\ref{lem:punctured-line}, they increase, and their union is $X$. If $X$ were Stein, so would be $M$, by the same lemma, contradicting~(a). Every $s\in H^0(M,K_M)$ defines a holomorphic function $F_s(x,\xi)=\langle\xi,s(x)\rangle$ on $X$. Points with different images $x\ne x'$ under $q$ are separated by $F_s$ for a section $s$ with $s(x)\ne0$ and $s(x')=0$, which exists by the Stein case of Theorem~\ref{thm:main}, and distinct points of one fiber of $q$ are separated by $F_s$ for any $s$ not vanishing at the base point. Choosing $s_0$ with $s_0(x)\ne0$ and $s_1,s_2$ such that $s_1/s_0,s_2/s_0$ are local coordinates at $x$ (the Stein case of Theorem~\ref{thm:main} with $k=1$), the functions $F_{s_0},F_{s_1},F_{s_2}$ give local coordinates at every point of $q^{-1}(x)$.
\end{proof}

Siu asked whether a complex manifold on which holomorphic functions and vector fields have arbitrary jets at any finite set is Stein if it satisfies the Kontinuit\"atssatz, or only its version for one-parameter families of disks~\cite[Question~(8.5)]{Si}. The space $(K_M^{-1})^\times$ answers the latter question negatively:

\begin{proposition}\label{prop:siu}
Holomorphic functions and holomorphic vector fields on $X=(K_M^{-1})^\times$ have arbitrary jets at any finite set. If $\varphi_t\colon\overline\Delta\to X$, $t\in[0,1)$, are holomorphic disks, continuous on $\overline\Delta$ and depending continuously on $t$, and $\bigcup_t\varphi_t(\partial\Delta)$ is relatively compact in $X$, then so is $\bigcup_t\varphi_t(\overline\Delta)$.
\end{proposition}

\begin{proof}
Let $P\subset X$ be finite. By the proof of Theorem~\ref{thm:fornaess}(b), finitely many functions in $\mathcal{O}(X)$ define a holomorphic map $G\colon X\to\mathbb{C}^N$ that is injective on $P$ and immersive at its points. Every germ at $p\in P$ is the pullback by $G$ of a germ at $G(p)$, and polynomials have arbitrary jets at the finite set $G(P)$, so the functions $Q\circ G$, $Q$ a polynomial, have arbitrary jets at $P$. In local coordinates $z$ of $M$ and the fiber coordinate $t$ of $K_M^{-1}$ with respect to the frame $\partial_{z_1}\wedge\partial_{z_2}$, a change of coordinates with Jacobian determinant $J$ replaces $t$ by $Jt$ and $dz_1\wedge dz_2\wedge dt$ by $J^2\,dz_1\wedge dz_2\wedge dt$. Hence $\Omega=t^{-2}dz_1\wedge dz_2\wedge dt$ is a nowhere vanishing holomorphic $3$-form on $X$, and $f,g\in\mathcal{O}(X)$ define a holomorphic vector field $V_{f,g}$ by $\iota_{V_{f,g}}\Omega=df\wedge dg$. If $f_1,f_2,f_3\in\mathcal{O}(X)$ are coordinates at $p$, then $V_{f_2,f_3},V_{f_3,f_1},V_{f_1,f_2}$ form a frame near $p$, and multiplying such frames by functions with prescribed jets at $P$ gives vector fields with arbitrary jets at $P$.

Let $X=\bigcup_jX_j$ with $X_j$ increasing and Stein, and let $K$ be the closure of $\bigcup_t\varphi_t(\partial\Delta)$. Choose $j$ with $K\cup\varphi_0(\overline\Delta)\subset X_j$. The set of $t$ such that $\varphi_s(\overline\Delta)\subset X_j$ for all $s\le t$ is open in $[0,1)$, and for such $t$ the maximum principle gives $\varphi_t(\overline\Delta)\subset\widehat K_{\mathcal{O}(X_j)}$, which is compact. By continuity, this set is also closed, so $\bigcup_t\varphi_t(\overline\Delta)\subset\widehat K_{\mathcal{O}(X_j)}$.
\end{proof}

\section{The C\oe ur\'e--Loeb space}\label{sec:CLsection}

\subsection{Graphs and suspensions}\label{sec:graph}
If $\phi\colon X\to N$ is locally Stein and $h\colon X\to Z$ is holomorphic, then $(\phi,h)\colon X\to N\times Z$ is locally Stein, since $(\phi,h)^{-1}(U\times Z)=\phi^{-1}(U)$ for every open set $U\subset N$.

Let $R$ be a Stein manifold and $T\in\operatorname{Aut}(R)$. The \emph{suspension} $E_T=(\mathbb{C}\times R)/\mathbb{Z}$, where $1\cdot(s,u)=(s+1,Tu)$, is a fiber bundle over $\mathbb{C}^*$ with projection $p[s,u]=e^{2\pi is}$. It carries the $\mathbb{C}$-action $\tau_v[s,u]=[s+v,u]$, which satisfies $p(\tau_ve)=e^{2\pi iv}p(e)$. If $B\subset\mathbb{C}^*$ is a disk and $\log$ is a branch of the logarithm on $B$, then $(b,u)\mapsto[(2\pi i)^{-1}\log b,u]$ is a biholomorphism $B\times R\to p^{-1}(B)$; hence $p$ is locally Stein.

\subsection{Cusps and C\oe ur\'e--Loeb bundles}\label{sec:cusp}

Fix $A=(a_{jk})=\bigl(\begin{smallmatrix}2&1\\1&1\end{smallmatrix}\bigr)\in\mathrm{SL}_2(\mathbb{Z})$, with eigenvalues $\lambda^{\pm1}$, where $\lambda>1$. Let $\ell_1$ be a left eigenvector for $\lambda$ with entries in $\mathbb{Q}(\lambda)$, and let $\ell_2$ be its Galois conjugate. We regard both as linear forms on $\mathbb{R}^2$. Put $C=\{\ell_1>0,\ell_2>0\}$; after multiplying $\ell_1$ by a suitable element of $\mathbb{Q}(\lambda)$, we can take $C$ to be any prescribed component of $\{\ell_1\ell_2\ne0\}$. The function $\operatorname{Nm}=\ell_1\ell_2$ is $A$-invariant on $C$, and $A$ acts properly discontinuously on $C$ with $\operatorname{Nm}\colon C/\langle A\rangle\to(0,\infty)$ proper. Let
\[
L(z)=(-\log|z_1|,-\log|z_2|),\qquad g_A(z)=(z_1^{a_{11}}z_2^{a_{12}},z_1^{a_{21}}z_2^{a_{22}}),
\]
so that $L\circ g_A=A\circ L$. The Reinhardt domain $D_C=L^{-1}(C)$ is Stein and biholomorphic to a bounded domain~\cite[p.~456]{OZ}. Its suspension $E=E_{A,C}$ is the \emph{C\oe ur\'e--Loeb space}; it is not Stein~\cite{CL}. The matrix used in~\cite{Pre26,Zaf01} is conjugate to $A$ by the interchange of the two coordinates.

\begin{lemma}\label{lem:cusp}
Let $Y_A=(D_C/\langle g_A\rangle)\cup\{\infty\}$, with the sets $U_c=(L^{-1}\{\operatorname{Nm}>c\}/\langle g_A\rangle)\cup\{\infty\}$, $c>0$, as neighborhoods of $\infty$. Then $Y_A$ is a normal Stein surface whose only singular point is the cusp $\infty$.
\end{lemma}

\begin{proof}
Write $z_j=e^{2\pi iw_j}$. Then $D_C/\langle g_A\rangle$ is the quotient of the tube $\{\operatorname{Im} w\in C\}$ by $\mathbb{Z}^2\rtimes\langle A\rangle$, and $(\ell_1,\ell_2)$ identifies this tube with $\mathbb{H}^2$. Under this identification, $\mathbb{Z}^2$ becomes a lattice $\mathfrak m\subset\mathbb{Q}(\lambda)$ acting by $(u_1,u_2)\mapsto(u_1+m,u_2+m')$, $A$ becomes multiplication by the totally positive unit $\lambda$, and $\operatorname{Nm}\circ L=(2\pi)^2\operatorname{Im} u_1\operatorname{Im} u_2$. By Hirzebruch~\cite[Section~2]{Hir73}, $Y_A$ is therefore a normal surface, and $\infty$ is a cusp singularity with neighborhood basis $U_c$.

The function $\psi=-\log h_1-\log h_2$, where the functions $h_j=\ell_j\circ L=2\pi\operatorname{Im}\ell_j(w)>0$ are pluriharmonic, satisfies $i\partial\bar\partial\psi=\sum_ji\,\partial h_j\wedge\bar\partial h_j/h_j^2>0$. It descends to $Y_A\setminus\{\infty\}$ and tends to $-\infty$ at $\infty$. The sets $\overline{U_c}=(L^{-1}\{\operatorname{Nm}\ge c\}/\langle g_A\rangle)\cup\{\infty\}$ are compact: for $c'>c$, the portion where $c\le\operatorname{Nm}\le c'$ is compact because $\operatorname{Nm}$ is proper on $C/\langle A\rangle$ and $L$ has compact fibers, while a sequence with $\operatorname{Nm}\to\infty$ converges to the added cusp point. Choose $c_1$ sufficiently large that a neighborhood of $\overline{U_{c_1}}$ admits an embedding $\iota$ into some $\mathbb{C}^m$, and choose $k>\max_{\overline{U_{c_1}}}|\iota|^2+\log(2c_1)$. Then
\[
\varphi=\begin{cases}
\max(\psi+k,|\iota|^2)&\text{on }U_{c_1},\\
\psi+k&\text{on }Y_A\setminus U_{c_1},
\end{cases}
\]
with $\varphi(\infty)=|\iota(\infty)|^2$, is continuous and strictly plurisubharmonic. Indeed, the first expression equals $\psi+k$ where $c_1<\operatorname{Nm}<2c_1$, and equals $|\iota|^2$ near $\infty$. Its sublevel sets $\{\varphi\le t\}\subset\overline{U_{e^{k-t}}}$ are compact, so $Y_A$ is Stein by Narasimhan's theorem~\cite{N62}.
\end{proof}

\begin{proof}[Proof of Theorem~\textup{\ref{thm:CL}}]
Let $E=E_{A,C}$ and regard $q\colon E\to Y_A$, $[s,u]\mapsto[u]$, as a map into the surface of Lemma~\ref{lem:cusp}. The map $\Pi=(p,q)\colon E\to\mathbb{C}^*\times Y_A$ is locally Stein by Section~\ref{sec:graph}, and its target is Stein. Its image is $\mathbb{C}^*\times Y_{A,\mathrm{reg}}$. To identify the covering, consider $\mathbb{C}\times D_C\to\mathbb{C}^*\times Y_{A,\mathrm{reg}}$, whose group $\mathbb{Z}^2$ acts by $(m,n)\cdot(s,u)=(s+m,g_A^nu)$. The map $\Pi$ is its quotient by the diagonal subgroup, so it is a regular covering of the regular locus with group $\mathbb{Z}^2/\mathbb{Z}(1,1)\cong\mathbb{Z}$. In particular, it is unramified.
\end{proof}

\begin{proposition}\label{prop:CL-canonical}
The canonical bundle of $E$ is holomorphically trivial.
\end{proposition}

\begin{proof}
Since $\det A=1$, the nowhere vanishing holomorphic $3$-form
\[
ds\wedge\frac{dz_1}{z_1}\wedge\frac{dz_2}{z_2}
\]
on $\mathbb{C}\times D_C$ is invariant under $(s,z)\mapsto(s+1,g_A(z))$. It therefore descends to a nowhere vanishing holomorphic $3$-form on $E$.
\end{proof}

\section{Locally Stein open subsets of singular Stein spaces}\label{sec:open}

\subsection{Thickening the graph}

\begin{proposition}\label{prop:strip}
Let $j\colon E_T\to Z$ be an open embedding into a complex space, and put $V=\{|\operatorname{Re} v|<1/4\}$. Then $F(e,v)=(p(e),j(\tau_ve))$ is an open embedding of $E_T\times V$ onto an open subset $W\subset\mathbb{C}^*\times Z_{\mathrm{reg}}$ that is locally Stein in $\mathbb{C}^*\times Z$. If $E_T$ is not Stein, then neither is $W$.
\end{proposition}

\begin{proof}
Suppose that $F(e,v)=F(e',v')$. Then $p(e)=p(e')$ and $\tau_ve=\tau_{v'}e'$. Applying $p$ gives $e^{2\pi iv}=e^{2\pi iv'}$; since $v\mapsto e^{2\pi iv}$ is injective on $V$, it follows that $v=v'$, and then $e=e'$. Thus $F$ is an injective holomorphic map between manifolds of the same dimension, hence an open embedding. The map $p\circ\mathrm{pr}_1\colon E_T\times V\to\mathbb{C}^*$ is locally Stein because $V$ is Stein, so $W$ is locally Stein by Section~\ref{sec:graph}. Finally, $F(E_T\times\{0\})\cong E_T$ is a closed submanifold of $W$, so $W$ is not Stein if $E_T$ is not.
\end{proof}

\subsection{An embedding with one singular point}\label{sec:embedding}

The following refinement of~\cite[Theorems~3.1 and~3.2]{Pre26} gives a normal Stein threefold with at most one singular point.

\begin{proposition}\label{prop:embed}
Let $Y$ be a connected normal Stein surface with one singular point $a$, let $\rho\colon\tilde Y\to Y$ be a resolution that is biholomorphic over $Y^*=Y\setminus\{a\}$, and let $P\to Y^*$ be a principal $(\mathbb{C},+)$-bundle. Suppose that its class $\xi\in H^1(Y^*,\mathcal{O})$ is the restriction of a class $\tilde\xi\in H^1(\tilde Y,\mathcal{O})$. Then $P$ is biholomorphic to the complement of an analytic subset of a normal Stein threefold $Z$ with at most one singular point.
\end{proposition}

\begin{proof}
We may assume that the exceptional curve $\mathcal C=\rho^{-1}(a)$ has simple normal crossings. It is connected since $Y$ is normal, and its intersection matrix $M$ is negative definite~\cite[p.~367]{Gra62}. Since $-M$ is positive definite with nonpositive off-diagonal entries, $(-M)^{-1}$ has positive entries, so a suitable positive integral multiple $D$ of $\sum_\Gamma((-M)^{-1}\mathbf 1)_\Gamma\,\Gamma$, where $\mathbf 1=(1,\dots,1)$, is an effective divisor supported on $\mathcal C$ with $D\cdot\Gamma<0$ for every irreducible component $\Gamma$ of $\mathcal C$. Then $\mathcal{O}(-D)$ has positive degree on every component of $\mathcal C$, so it is ample on $\mathcal C$, hence on the fiber $\rho^{-1}(a)$, whose reduction is $\mathcal C$~\cite[Corollary~IV.3.2 and Exercise~III.5.7]{Har77}, and thus on every fiber of $\rho$, the other fibers being points. After replacing $D$ by a multiple, the relative Theorems~A and~B for fiberwise ample line bundles~\cite[Theorem~3.2.4]{Con06} give
\[
\rho^*\rho_*\mathcal{O}(-D)\longrightarrow\mathcal{O}(-D)\ \text{surjective},\qquad R^1\rho_*\mathcal{O}(-D)=0
\]
near $a$. Conrad's proof applies verbatim to complex spaces~\cite[Introduction]{Con06}; for locally projective maps, see Grauert and Remmert~\cite{GR58}. Away from $a$, these assertions hold because $\rho$ is biholomorphic there. Cartan's Theorem~B and the Leray spectral sequence give $H^1(\tilde Y,\mathcal{O}(-D))=0$.

The sheaf $\mathcal{J}=\rho_*\mathcal{O}(-D)\subset\mathcal{O}_Y$ is a coherent ideal, equal to $\mathcal{O}_Y$ off $a$. Choose finitely many $f_1,\dots,f_r\in\Gamma(Y,\mathcal{J})$ that generate $\mathcal{J}_a$ and have no common zero on $Y^*$. To achieve the latter condition, start with global sections generating $\mathcal{J}_a$ (Cartan's Theorem~A), and let $\Sigma$ be their common zero set. Since they generate $\mathcal{J}$ near $a$ and $\mathcal{J}=\mathcal{O}_Y$ off $a$, the point $a$ is isolated in $\Sigma$, so $\Sigma'=\Sigma\setminus\{a\}$ is analytic and $a\notin \Sigma'$. The sheaf $\mathcal{J}/\mathcal{J}\mathcal I_{\Sigma'}$ equals $\mathcal{O}_{\Sigma'}$, because $\mathcal{J}=\mathcal{O}_Y$ near $\Sigma'$. By Cartan's Theorem~B, the section $1$ of $\mathcal{O}_{\Sigma'}$ lifts to some $f\in\Gamma(Y,\mathcal{J})$, and adjoining $f$ removes all common zeros on $Y^*$. The pullbacks $\rho^*f_\nu$ therefore generate $\mathcal{O}(-D)$ everywhere.

Let $\{U_i\}$ be a Stein cover of $\tilde Y$ on which $D$ has local equations $h_i$, and represent $\tilde\xi$ by a cocycle $\tilde c_{ij}$. Write $\rho^*f_\nu=h_i\phi_{\nu,i}$. Since $H^1(\tilde Y,\mathcal{O}(-D))=0$, there are sections $b_{\nu,i}\in\Gamma(U_i,\mathcal{O}(-D))$, viewed as holomorphic functions vanishing on $D$, such that
\[
b_{\nu,i}-b_{\nu,j}=\rho^*f_\nu\tilde c_{ij}.
\]
Define an affine line bundle $P'\to\tilde Y$ by fiber coordinates $t_i$ with transition maps
\[
t_i=(h_i/h_j)t_j+h_i\tilde c_{ij}.
\]
Over $\tilde Y\setminus\mathcal C$, the coordinates $s_i=t_i/h_i$ identify $P'$ with $P$. Over $\mathcal C$, the translations vanish, and $P'|_{\mathcal C}$ is the total space of $\mathcal{O}(D)|_{\mathcal C}$. The functions
\[
F_\nu=\phi_{\nu,i}t_i-b_{\nu,i}
\]
glue to functions on $P'$. The map $\Phi'=(\rho,F_1,\dots,F_r)\colon P'\to Y\times\mathbb{C}^r$ is proper: over the inverse image of a compact subset of $Y$, finitely many local trivializations suffice, and since the $\phi_\nu$ have no common zero, the fiber coordinates are bounded whenever all $F_\nu$ are bounded. Over $Y^*$ the map is injective, since some $f_\nu$ is nonzero at each point.

On $P'|_{\mathcal C}$ we have $F_\nu=\phi_\nu t$. The zero section $\mathcal C_0$ is precisely the fiber over $(a,0)$, and every other fiber is finite. Indeed, the projective map $y\mapsto[\phi_1(y):\cdots:\phi_r(y)]$ is nonconstant on each component $\Gamma$ of $\mathcal C$: if it were constant, the generating sections of $\mathcal{O}(-D)|_\Gamma$, a line bundle of positive degree, would be multiples of a single nowhere vanishing section, which is impossible. Each fiber of this projective map is therefore finite, and a nonzero value of $F$ then determines $t$ uniquely.

Let $P'\xrightarrow{\psi}Z\to Y\times\mathbb{C}^r$ be the Stein factorization of $\Phi'$~\cite{GR84}: the space $Z$ is normal and finite over the Stein target, hence Stein, $\psi_*\mathcal{O}_{P'}=\mathcal{O}_Z$, and the fibers of $\psi$ are the connected components of those of $\Phi'$. Hence $\psi$ contracts the connected curve $\mathcal C_0$ to a point and is injective elsewhere, so it is biholomorphic off $\mathcal C_0$. Consequently, $Z_{\mathrm{sing}}\subset\psi(\mathcal C_0)$ and $P\cong Z\setminus\psi(P'|_{\mathcal C})$, where the removed set is analytic since $\psi$ is proper. Since $P'$ is a connected three-dimensional manifold and $\psi$ is biholomorphic off $\mathcal C_0$, $Z$ is a threefold.
\end{proof}

\begin{theorem}\label{thm:CLembed}
The C\oe ur\'e--Loeb space $E$ is biholomorphic to an open subset of a normal Stein threefold with exactly one singular point.
\end{theorem}

\begin{proof}
The map $q\colon E\to Y_A\setminus\{\infty\}$, $[s,u]\mapsto[u]$, is a principal $(\mathbb{C},+)$-bundle for the action $\tau$. Zaffran extends this bundle across the exceptional cycle of a resolution of the cusp~\cite[Section~7]{Zaf01}; this is the geometric input of Preda's theorem~\cite[Section~2.2 and Theorem~3.2]{Pre26}. Proposition~\ref{prop:embed} therefore gives an open embedding $E\subset Z$ into a normal Stein threefold $Z$ with at most one singular point. If $Z$ were smooth, Proposition~\ref{prop:strip} would give a locally Stein, non-Stein open subset of the Stein manifold $\mathbb{C}^*\times Z$, contradicting~\cite{DG}. Thus $Z$ has exactly one singular point.
\end{proof}

\subsection{A compact set in no Stein open subset}

\begin{proposition}\label{prop:disc}
The C\oe ur\'e--Loeb space $E$ contains a compact connected set that lies in no Stein open subset.
\end{proposition}

\begin{proof}
Following C\oe ur\'e and Loeb~\cite{CL}, we construct disks in the coordinates of Lemma~\ref{lem:cusp}. The space $E$ is covered by $\mathbb{C}\times\mathbb{H}^2$, with deck transformations given by real lattice translations in $u=(u_1,u_2)$ and
\[
(s,u_1,u_2)\longmapsto(s+1,\lambda u_1,\lambda^{-1}u_2).
\]
Put $\mu=\log\lambda$. For $1<R\le2$ define
\[
f_R(\zeta)=\mu^{-1}\log\frac{R+\zeta}{R-\zeta}\qquad(\zeta\in\overline\Delta),
\]
using the logarithm on the right half-plane. Let $u_{1,R}$ and $u_{2,R}$ be holomorphic on $\Delta$ and continuous on $\overline\Delta$, with boundary imaginary parts $e^{\mu\operatorname{Re}f_R}$ and $e^{-\mu\operatorname{Re}f_R}$, respectively, and normalized by $\operatorname{Re}u_{j,R}(0)=0$. They are obtained by harmonic extension and harmonic conjugation, take values in $\mathbb{H}$, and depend continuously on $R$. The maps $(f_R,u_{1,R},u_{2,R})$ thus project to a continuous family of closed holomorphic disks $d_R$ in $E$.

Their boundaries lie in the compact connected set $Q_0\subset E$ given by the image of the convex set of points $(s,u)$ with $0\le\operatorname{Re}s\le1$, $|\operatorname{Im}s|\le\pi/(2\mu)$, $1\le\operatorname{Im}u_1\le\lambda$, $\lambda^{-1}\le\operatorname{Im}u_2\le1$ and $(\operatorname{Re}u_1,\operatorname{Re}u_2)$ in a fixed compact fundamental parallelogram of the lattice translations. Indeed, on the boundary $\operatorname{Im}u_1=e^{\mu\operatorname{Re}s}$ and $\operatorname{Im}u_2=e^{-\mu\operatorname{Re}s}$, and an integral power of the displayed deck transformation, which preserves these relations, puts $\operatorname{Re}s$ in $[0,1]$; moreover $|\operatorname{Im}s|\le\pi/(2\mu)$, and lattice translations adjust $\operatorname{Re}u$. At the center, $f_R(0)=0$ and $u_{1,R}(0)=u_{2,R}(0)=ia_R$, where
\[
a_R=\frac1{2\pi}\int_0^{2\pi}\left|\frac{R+e^{it}}{R-e^{it}}\right|\,dt\longrightarrow\infty\quad(R\downarrow1).
\]
The equality of the two means follows by replacing $t$ by $t+\pi$, and divergence follows from Fatou's lemma. The invariant function $\operatorname{Im}u_1\operatorname{Im}u_2$ therefore tends to infinity at $d_R(0)$, so the centers leave every compact subset of $E$.

Put $Q=Q_0\cup d_2(\overline\Delta)$, a compact connected set, and suppose that $Q\subset U$ for some Stein open set $U$. The set $I=\{R:d_R(\overline\Delta)\subset U\}$ contains $2$ and is open by continuity. By the maximum principle, every disk with $R\in I$ lies in the compact hull $\widehat Q_{\mathcal{O}(U)}$. At a limit parameter in $(1,2]$, continuity keeps the whole disk in this hull, so $I$ is also closed. Hence $I=(1,2]$, contradicting the escape of the centers.
\end{proof}

\subsection{Proof of Theorem~\ref{thm:open}}

\begin{proof}[Proof of Theorem~\ref{thm:open}]
Let $Z$ be the normal Stein threefold with one singular point from Theorem~\ref{thm:CLembed}, which contains the C\oe ur\'e--Loeb space $E$ as an open subset, and put $N=\mathbb{C}^*\times Z$. By Proposition~\ref{prop:strip}, $W=F(E\times V)$ is a locally Stein, non-Stein open subset of $N_{\mathrm{reg}}$.

For the relatively compact version, let $Q\subset E$ be as in Proposition~\ref{prop:disc}. Then $Q'=F(Q\times\{0\})$ is compact and connected, and lies in no Stein open subset of $W$: intersecting such an open set with the closed complex submanifold $F(E\times\{0\})$ would give a Stein open subset of $E$ containing $Q$.

Let $\varphi$ be a continuous strictly plurisubharmonic exhaustion function of $N$, and let $c>\max_{Q'}\varphi$. The set $N_{<c}=\{\varphi<c\}$ is Stein and relatively compact in $N$. Let $W'$ be the connected component of $W\cap N_{<c}$ containing $Q'$; it is relatively compact and not Stein.

For each point of $N$, choose an open neighborhood $U$ such that $U\cap W$ is Stein. The intersection $U\cap W\cap N_{<c}$ is Stein, since the diagonal of $N$ identifies it with a closed analytic subspace of $(U\cap W)\times N_{<c}$. The set $U\cap W'$ is a union of connected components of $U\cap W\cap N_{<c}$, and hence is Stein. Thus $W'$ is locally Stein. Since $W'\subset W\subset N_{\mathrm{reg}}$, the set $W'$ has all the properties required of $W$ in Theorem~\ref{thm:open}.
\end{proof}

\begin{remark}\label{rem:gluing}
Col\c{t}oiu observed that the Levi problem for open subsets of Stein spaces is equivalent to asking whether an open subset $D$ of a Stein space $Y$ is Stein if, for some $f\in\mathcal{O}(Y)$ and $a<b$, the sets $D\cap\{\operatorname{Re}f<b\}$ and $D\cap\{\operatorname{Re}f>a\}$ are Stein~\cite[Remark~6.2]{C09}. The set $W=F(E\times V)$ in the proof above is a counterexample, with $Y=N$, $f(t,z)=t+t^{-1}$, $a=-1$ and $b=1$. Indeed, $S_-=\{t\in\mathbb{C}^*:\operatorname{Re}(t+t^{-1})<1\}$ avoids the positive real axis, and $S_+=\{t\in\mathbb{C}^*:\operatorname{Re}(t+t^{-1})>-1\}$ avoids the negative one. Hence a branch of the logarithm trivializes $p$ over each of them, and since the first component of $F$ is $p$,
\[
W\cap(S_\pm\times Z)=F\bigl(p^{-1}(S_\pm)\times V\bigr)\cong S_\pm\times D_C\times V,
\]
which is Stein.
\end{remark}

\subsection*{Acknowledgments}
The first author was supported by JSPS KAKENHI Grant Number JP24K16919.

\subsection*{AI disclosure}
The main mathematical ideas and proof strategies originated with the authors. Following these strategies, most of the detailed proofs in Sections~\ref{sec:prelim}--\ref{sec:open} were worked out with the assistance of generative AI tools. These tools were also used to check arguments, to search the literature, and to prepare and revise the manuscript. The authors take full responsibility for the content of this paper.

\end{document}